\documentclass[12pt]{article}
\usepackage{amssymb,amsmath,amsthm}
\renewcommand*{\leq}{\leqslant} \renewcommand*{\geq}{\geqslant}

\newcommand*{\vare}{\varepsilon} \newcommand*{\vart}{\vartheta}

\newcommand*{\mC}{\mathbb C} \newcommand*{\mN}{\mathbb N}
\newcommand*{\mR}{\mathbb R} \newcommand*{\mS}{\mathbb S}
\newcommand*{\mT}{\mathbb T} \newcommand*{\mZ}{\mathbb Z}

\newcommand*{\cA}{\mathcal A} \newcommand*{\cC}{\mathcal C}
\newcommand*{\cF}{\mathcal F} \newcommand*{\cG}{\mathcal G}
\newcommand*{\cL}{\mathcal L} \newcommand*{\cM}{\mathcal M}
\newcommand*{\cN}{\mathcal N} \newcommand*{\cO}{\mathcal O}
\newcommand*{\cT}{\mathcal T} \newcommand*{\cX}{\mathcal X}

\newcommand*{\fc}{\mathfrak c} \newcommand*{\fg}{\mathfrak g}
\newcommand*{\fl}{\mathfrak l} \newcommand*{\fp}{\mathfrak p}
\newcommand*{\fs}{\mathfrak s} \newcommand*{\ft}{\mathfrak t}
\newcommand*{\fu}{\mathfrak u} \newcommand*{\fw}{\mathfrak w}
\newcommand*{\fx}{\mathfrak x} \newcommand*{\fy}{\mathfrak y}
\newcommand*{\fz}{\mathfrak z}

\newcommand*{\fA}{\mathfrak A} \newcommand*{\fD}{\mathfrak D}
\newcommand*{\fG}{\mathfrak G} \newcommand*{\fM}{\mathfrak M}
\newcommand*{\fO}{\mathfrak O} \newcommand*{\fP}{\mathfrak P}
\newcommand*{\fS}{\mathfrak S} \newcommand*{\fT}{\mathfrak T}
\newcommand*{\fV}{\mathfrak V} \newcommand*{\fW}{\mathfrak W}
\newcommand*{\fX}{\mathfrak X} \newcommand*{\fY}{\mathfrak Y}
\newcommand*{\fZ}{\mathfrak Z}

\newcommand*{\bx}{\overline{x}} \newcommand*{\by}{\overline{y}}
\newcommand*{\bz}{\overline{z}}

\newcommand*{\gl}{\mathfrak{gl}}

\newcommand*{\GL}{\mathrm{GL}} \newcommand*{\SL}{\mathrm{SL}}

\newcommand*{\hM}{\widehat{M}} \newcommand*{\hX}{\widehat{X}}
\newcommand*{\hY}{\widehat{Y}} \newcommand*{\hZ}{\widehat{Z}}

\newcommand*{\halpha}{\widehat{\alpha}} \newcommand*{\hbeta}{\widehat{\beta}}

\newcommand*{\hTheta}{\widehat{\Theta}} \newcommand*{\hOmega}{\widehat{\Omega}}

\newcommand*{\tA}{\widetilde{A}} \newcommand*{\tB}{\widetilde{B}}
\newcommand*{\tM}{\widetilde{M}}

\newcommand*{\talpha}{\widetilde{\alpha}}
\newcommand*{\tbeta}{\widetilde{\beta}}

\newcommand*{\tGamma}{\widetilde{\Gamma}}
\newcommand*{\tOmega}{\widetilde{\Omega}}

\newcommand*{\bfzero}{\mathbf 0}

\newcommand*{\bfomega}{\text{\boldmath $\omega^2$}}

\newcommand*{\rme}{\mathrm e} \newcommand*{\rmi}{\mathrm i}

\newcommand*{\add}{^{\mathrm{add}}} \newcommand*{\new}{^{\mathrm{new}}}

\DeclareMathOperator{\codim}{codim} \DeclareMathOperator{\meas}{meas}

\DeclareMathOperator{\Fix}{Fix}

\newtheorem{lmm}{Lemma} \newtheorem{prp}{Proposition}
\newtheorem{thm}{Theorem}

\theoremstyle{definition}
\newtheorem{dfn}{Definition} \newtheorem{xpl}{Example}
\newtheorem{rmk}{Remark}

\allowdisplaybreaks

\begin{document}
\baselineskip=21pt
\parindent=0mm

\begin{center}
{\Large\bfseries Partial Preservation \\
of Frequencies and Floquet Exponents of Invariant Tori \\
in the Reversible KAM Context~2}

\bigskip
\baselineskip=16pt

{\large\bfseries Mikhail B.~Sevryuk\footnote{E-mails: \texttt{sevryuk@mccme.ru, 2421584@mail.ru}}}

\medskip

{\small\itshape V.~L.~Talroze Institute of Energy Problems of Chemical Physics of the Russia Academy of Sciences, \\
Leninskii prospect~38, Building~2, Moscow 119334, Russia}
\end{center}

\bigskip
\baselineskip=17pt

\textbf{Abstract}---%
We consider the persistence of smooth families of invariant tori in the reversible context~2 of KAM theory under various weak nondegeneracy conditions via Herman's method. The reversible KAM context~2 refers to the situation where the dimension of the fixed point manifold of the reversing involution is less than half the codimension of the invariant torus in question. The nondegeneracy conditions we employ ensure the preservation of any prescribed subsets of the frequencies of the unperturbed tori and of their Floquet exponents (the eigenvalues of the coefficient matrix of the variational equation along the torus).

\bigskip

MSC2010 numbers: \texttt{70K43, 70H33}

\bigskip


\bigskip

\begin{flushright}\itshape
To the blessed memory of Helmut R\"ussmann \textup{(}1930--2011\textup{)} \\
whose contribution to KAM theory is so substantial and versatile
\end{flushright}

\parindent=5mm

\section{Introduction}\label{introduction}

\subsection{Reversible Contexts~1 and~2}\label{reversiblecontexts}

Equilibria, periodic orbits, invariant tori filled up with quasi-periodic motions (conditionally periodic motions with rationally independent frequencies) and their asymptotic manifolds (in particular, homoclinic and heteroclinic trajectories) are key elements of finite dimensional dynamics. The importance of equilibria (invariant $0$-tori) and periodic orbits (invariant $1$-tori) of autonomous flows was realized already by H.~Poincar\'e and further emphasized from the bifurcational viewpoint by A.A.~Andronov and E.~Hopf \cite{MM76}. Quasi-periodic motions with $n\geq 2$ basic frequencies are the subject of the Kolmogorov--Arnold--Moser (KAM) theory founded in the fifties and sixties of the last century \cite{AKN06,BHS96LNM,BS10,dL01,D14,GHdL14,HCFLM16,M67MA,R01RCD,XL16RCD}. According to KAM theory, the occurrence of invariant tori of various dimensions carrying quasi-periodic motions and organized into Cantor-like families is a generic property of non-integrable dynamical systems. The possible dimensions of the tori and the number of parameters of their Cantor families (as a rule, these families themselves form complicated hierarchical conglomerates) depend strongly on the phase space structures the system in question is assumed to preserve.

For instance, a typical autonomous Hamiltonian system with $N$ degrees of freedom is expected to admit isolated equilibria, one-parameter smooth families of periodic orbits (the parameter being just the energy value), and $n$-parameter Cantor families of isotropic invariant $n$-tori carrying quasi-periodic motions for each $n=2,\ldots,N$ \cite{AKN06,BHS96LNM,BS10,R01RCD}. The existence of other types of families of quasi-periodic motions filling up isotropic invariant tori is an evidence for the presence of additional symmetries of the system. By the way, in a generic one-parameter family of periodic orbits, the period is not a constant and can be used as an alternative parameter.

In KAM theory, one considers various classes of dynamical systems, and the invariant tori sought for can relate to the corresponding phase space structures in different ways, so one sometimes speaks of particular \emph{contexts} of KAM theory. The most explored finite dimensional contexts are the dissipative context (with no special structures on the phase space), the volume preserving context (where one looks for invariant tori of volume preserving systems), the Hamiltonian isotropic context (where one examines isotropic invariant tori in Hamiltonian systems), and the so-called reversible context~1 \cite{BHN07JDE,BH95JDDE,BHS96G,BHS96LNM,BHT90,QS93C,S95UMN,S95C,S06N,S07DCDS,S07TMIS,W10DCDSS}. Less familiar contexts are exemplified by the Hamiltonian coisotropic context (with coisotropic invariant tori) and the Hamiltonian atropic context (where the invariant tori to be constructed are atropic, i.e., neither isotropic nor coisotropic), see \cite{BS10} for references on the Hamiltonian coisotropic and atropic contexts, as well as by the so-called conformally Hamiltonian context, see \cite{CCdL17N} and references therein. One more example is the reversible context~2 the present paper is devoted to. Let us recall the relevant definitions and principal facts.

\begin{dfn}[\cite{LR98PD,RQ92PR,S86}]\label{reversible}
Given an arbitrary set $\cM$, a mapping $G:\cM\to\cM$ is called an \emph{involution} of $\cM$ if $G^2=G\circ G$ is the identity transformation. A dynamical system is said to be \emph{reversible} with respect to a smooth involution $G$ of the phase space (or $G$-reversible) if this system is invariant under the transformation $(\fp,t)\mapsto(G\fp,-t)$ where $\fp$ is a point of the phase space and $t$ is the time (i.e., if $G$ casts the system in question into the system with the reverse time direction).
\end{dfn}

In the reversible KAM theory, one always deals with only those tori that are invariant under both the system itself and the reversing involution.

\begin{lmm}[\cite{BHS96G,BHS96LNM,S86,S12MMJ}]\label{torus}
Let an $n$-torus $\cT\subset\cM$ be invariant under both a $G$-reversible flow on $\cM$ and the corresponding reversing involution $G$. If $\cT$ carries quasi-periodic motions then one can introduce a coordinate frame $x\in\mT^n=(\mR/2\pi\mZ)^n$ in $\cT$ such that the dynamics on $\cT$ takes the form $\dot{x}=\omega$ and the restriction of $G$ to $\cT$ takes the form $G|_{\cT}:x\mapsto-x$. Consequently, the set of fixed points of $G|_{\cT}$ consists of $2^n$ isolated points $(x_1,\ldots,x_n)$ where each $x_i$, $1\leq i\leq n$, is equal either to $0$ or to $\pi$.
\end{lmm}

The set $\Fix G$ of fixed points of an involution $G:\cM\to\cM$ of a manifold $\cM$ is a submanifold of $\cM$ of the same smoothness class as $G$ itself \cite{B72,MZ74} (the books \cite{B72,CF64} present extensive information on the fixed point sets of involutions of various manifolds). However, in the framework of Lemma~\ref{torus}, different points of $\Fix(G|_{\cT})=(\Fix G)\cap\cT$ may belong to connected components of $\Fix G$ of different dimensions (see \cite{QS93C,S95UMN} for several examples in the case $n=1$). None of these dimensions can exceed the codimension $\codim\cT$ of $\cT$ in the phase space because $\dim(\cT\cap\Fix G)=0$.

\begin{dfn}[\cite{BHS96G,BHS96LNM}]\label{contexts}
Let all the connected components of $\Fix G$ that intersect $\cT$ in the framework of Lemma~\ref{torus} be of the same dimension $d_G$. The situation where the inequalities $\frac{1}{2}c_{\cT}\leq d_G\leq c_{\cT}$ hold (here $c_{\cT}=\codim\cT$) is called \emph{the reversible context~1}. The opposite situation where the inequality $d_G<\frac{1}{2}c_{\cT}$ holds is called \emph{the reversible context~2}.
\end{dfn}

Note that for most involutions $G$ encountered in practice, the fixed point manifold $\Fix G$ is not empty and all the connected components of $\Fix G$ are of the same dimension, so $\dim\Fix G$ is well defined \cite{LR98PD,RQ92PR}.

The drastic differences between the two reversible contexts and the peculiarities of the reversible context~2 are discussed in detail in our previous articles \cite{S11RCD,S12MMJ,S12IM,S16RCD,S17MMJ}. Here we just demonstrate these differences in the trivial case $n=0$ where the invariant tori in question are equilibria and their codimension is the phase space dimension. These equilibria should be invariant under the reversing involution $G$, i.e., they should be fixed points of $G$.

\begin{xpl}\label{trivial}
Consider the involution $G:(u,v)\mapsto(u,-v)$ of $\mR^{a+b}$ where $u\in\mR^a$ and $v\in\mR^b$, so that $\Fix G=\{v=0\}$ and $\dim\Fix G=a$. A system
\[
\dot{u}=U(u,v), \qquad \dot{v}=V(u,v)
\]
is reversible with respect to $G$ if and only if $U$ is odd in $v$ and $V$ is even in $v$. We are looking for the equilibria of such a system on the plane $\Fix G$, i.e., for the points $u\in\mR^a$ such that $U(u,0)=0$ and $V(u,0)=0$. However, $U(u,0)\equiv 0$, and the desired equilibria $(u,0)$ are determined by the equation $V(u,0)=0$.

The reversible context~1 here corresponds to the inequality $\frac{1}{2}(a+b)\leq a$, i.e., $a\geq b$. Within this context, the equation $V(u,0)=0$ generically describes a smooth $(a-b)$-dimensional surface in $\Fix G$. On the other hand, in the reversible context~2 (where $a<b$) one generically has no equilibria lying in $\Fix G$. To obtain such equilibria, one has to let the system depend on at least $b-a$ \emph{external parameters}. For a $G$-reversible system $\dot{u}=U(u,v,\fw)$, $\dot{v}=V(u,v,\fw)$ depending on a $\fc$-dimensional external parameter $\fw$ with $\fc\geq b-a$, one generically gets a smooth $(\fc-b+a)$-dimensional surface of equilibria in the product of the plane $\Fix G$ and the parameter space $\mR^{\fc}\ni\fw$.
\end{xpl}

Let $R\in\GL(a+b,\mR)$ be an involutive matrix with eigenvalue $1$ of multiplicity $a$ and eigenvalue $-1$ of multiplicity $b$. One says that a matrix $M\in\gl(a+b,\mR)$ anti-commutes with $R$, or is \emph{infinitesimally reversible} with respect to $R$, if $MR=-RM$. If this is the case then the eigenvalues of $M$ come in pairs $(\lambda,-\lambda)$, and if $b\neq a$ then $0$ is an eigenvalue of $M$ of multiplicity $\ft\geq|b-a|$ \cite{H96JDE,S86,S91TSP,S93CJM} (generically $\ft=|b-a|$).

Consider the linearization of a $G$-reversible system in the setup of Example~\ref{trivial} around any equilibrium lying in $\Fix G$. If $b\neq a$ then this linearization possesses the zero eigenvalue of multiplicity $\ft\geq|b-a|$ (generically $\ft=|b-a|$). The nonzero eigenvalues come in pairs $(\lambda,-\lambda)$.

\subsection{Unperturbed Systems in the Reversible Context~2}\label{setup}

It is an appropriate time now to introduce some notation. Let $\mN$ be the set of positive integers and let $\mZ_+=\mN\cup\{0\}$. Throughout the paper, we will denote by $|{\cdot}|$ the $\ell_1$-norm of vectors in $\mC^s$, by $\|{\cdot}\|$ the $\ell_2$-norm of vectors in $\mR^s$, and by $\langle{\cdot},{\cdot}\rangle$ the inner product of two vectors in $\mR^s$. A closed $s$-dimensional ball centered at a point $\mu\in\mR^s$ is the set $B=\bigl\{ \fp\in\mR^s \bigm| \|\fp-\mu\|\leq\varrho \bigr\}$ for a certain $\varrho>0$. For $s=0$ this definition gives $\mu=0$ and $B=\{0\}=\mR^0$. The expression $\cO_s(\mu)$ will denote an unspecified neighborhood of a point $\mu\in\mR^s$. If $d\in\mN$ and $x,y,z,\ldots$ are certain variables, we will write $O_d(x,y,z,\ldots)$ instead of $O\bigl( |x|^d+|y|^d+|z|^d+\cdots \bigr)$. Instead of $O_1({\cdot})$, we will write just $O({\cdot})$.

Given a matrix $M\in\gl(N,\mR)$, the expression $\bfzero_m\oplus M$ will denote the $(m+N)\times(m+N)$ block diagonal matrix whose first block is the $m\times m$ zero matrix and the second block is $M$. The space of $n\times N$ real matrices will be denoted by $\mR^{n\times N}$, so that $\gl(n,\mR)=\mR^{n\times n}$.

Recall also that a $C^1$-smooth mapping $\cF:\cM\to\cN$ of smooth manifolds is said to be \emph{submersive} at a point $\mu\in\cM$, if $\dim\cM\geq\dim\cN$ and the rank of the differential of $\cF$ is equal to $\dim\cN$ at $\mu$. If this is the case then $\cF$ is also submersive at any point $\mu'\in\cM$ sufficiently close to $\mu$.

\begin{dfn}\label{reducible}
Let $\cT$ be an invariant $n$-torus of some flow on an $(n+N)$-dimensional manifold. This torus is said to be \emph{reducible} (or \emph{Floquet}) if in a neighborhood of $\cT$, there exists a coordinate frame $x\in\mT^n$, $\cX\in\cO_N(0)$ in which the torus $\cT$ itself is given by the equation $\cX=0$ and the dynamical system takes the \emph{Floquet form} $\dot{x}=\omega+O(\cX)$, $\dot{\cX}=\Lambda\cX+O_2(\cX)$ with $x$-independent vector $\omega\in\mR^n$ and matrix $\Lambda\in\gl(N,\mR)$. The vector $\omega$ (not determined uniquely) is called the \emph{frequency vector} of the torus $\cT$, while the matrix $\Lambda$ (not determined uniquely) is called the \emph{Floquet matrix} of $\cT$, and its eigenvalues are called the \emph{Floquet exponents} of $\cT$. The coordinates $(x,\cX)$ are called the \emph{Floquet coordinates} for $\cT$.
\end{dfn}

Note that the Floquet exponents of an equilibrium (where $n=0$) are just the eigenvalues of the linearization of the vector field around this equilibrium.

In the overwhelming majority of works on KAM theory, the invariant tori under study are reducible. In particular, this is the case for all the papers on the reversible context~2 \cite{S11RCD,S12MMJ,S12IM,S16RCD,S17MMJ}. The Cantor families of reducible invariant tori in KAM theory are in fact \emph{Whitney smooth}. This means that although the Floquet coordinates for the tori within a given $\fs$-parameter family are defined a~priori on a certain Cantor-like subset of $\mR^{\fs}$, these coordinates can be continued to smooth (say, $C^\infty$) functions defined in an open domain in $\mR^{\fs}$. For basic references on Whitney smoothness in KAM theory, see \cite{BHS96LNM,BS10}.

The results of the present paper imply that the situation with reducible invariant tori of an arbitrary dimension $n$ within the reversible context~2 is more or less similar to the trivial case $n=0$ of Example~\ref{trivial}. Namely, if the phase space codimension of each torus is equal to $a+b$ and $\dim\Fix G=a<b$ ($G$ being the reversing involution), then for $n\geq 2$ one needs at least $b-a+1$ external parameters, to be more precise, $b-a$ parameters for the same reasons as in the case $n=0$ (cf.\ Proposition~\ref{notorus} in Section~\ref{multiparameter} below) and one more parameter to control resonances involving the frequencies and the imaginary parts of the Floquet exponents. Each torus possesses the zero Floquet exponent of multiplicity $\ft\geq b-a$ (generically $\ft=b-a$). The nonzero Floquet exponents come in pairs $(\lambda,-\lambda)$. If the number of external parameters is equal to $\fc\geq b-a+1$, one obtains a $(\fc-b+a)$-parameter Cantor family of invariant tori in the product of the phase space and the space of external parameters. It is worthwhile to emphasize that in the four ``conventional'' KAM contexts (the reversible context~1, Hamiltonian isotropic context, volume preserving context, and dissipative context), a one-dimensional external parameter is always enough \cite{BHS96G,BHS96LNM,S95C,S06N,S07DCDS,S07TMIS} (with the exception of very special situations).

In fact, to run KAM theory for the reversible context~2, one has first to choose the unperturbed systems where the invariant tori are organized into a $(\fc-b+a)$-parameter smooth (rather than Cantor) family. Following \cite{S16RCD,S17MMJ}, we will consider unperturbed systems of the form
\begin{equation}
\begin{aligned}
\dot{x} &= \Omega(\mu)+\Delta(\sigma,\mu)+\xi(y,z,\sigma,\mu), \\
\dot{y} &= \sigma+\eta(y,z,\sigma,\mu), \\
\dot{z} &= M(\mu)z+\zeta(y,z,\sigma,\mu),
\end{aligned}
\label{unperturbed}
\end{equation}
where $x\in\mT^n$, $y\in\cO_m(0)$, $z\in\cO_{2p}(0)$ are the phase space variables, $\sigma\in\cO_m(0)$ and $\mu\in\cO_s(0)$ are external parameters ($n\in\mZ_+$, $m\in\mN$, $p\in\mZ_+$, $s\in\mN$), $M$ is a $2p\times 2p$ matrix-valued function, $\Delta=O(\sigma)$ and $\xi=O(y,z)$, $\eta=O_2(y,z)$, $\zeta=O_2(y,z,\sigma)$. These systems are assumed to be reversible with respect to the involution
\begin{equation}
G:(x,y,z)\mapsto(-x,-y,Rz),
\label{involution}
\end{equation}
where $R\in\GL(2p,\mR)$ is an involutive matrix with eigenvalues $1$ and $-1$ of multiplicity $p$ each and $M(\mu)R\equiv-RM(\mu)$. The dimension of the space $\bigl\{(\sigma,\mu)\bigr\}$ of external parameters is equal to $m+s$. The systems~\eqref{unperturbed} are ``integrable'' in the sense that they are $\mT^n$-equivariant, i.e., the right-hand side of~\eqref{unperturbed} is independent of the angular variable $x$.

For $\sigma=0$ and any value of $\mu$, the system~\eqref{unperturbed} and involution~\eqref{involution} admit a common reducible invariant $n$-torus $\{y=0, \, z=0\}$ with frequency vector $\Omega(\mu)\in\mR^n$ and Floquet matrix $\bfzero_m\oplus M(\mu)\in\gl(m+2p,\mR)$. The codimension of this torus in the phase space is equal to $m+2p$ while $\dim\Fix G=p$ (in the previous notation, $a=p$, $b=m+p>a$, and $\fc=m+s$, so that $\fc-b+a=s$). The parameter $\sigma$ is a remedy for a drift along the variable $y$ in $G$-reversible perturbations of the systems~\eqref{unperturbed}.

\begin{rmk}\label{Pi}
One may wonder why the equation for $\dot{z}$ in~\eqref{unperturbed} does not contain a term like $\Pi(\sigma,\mu)z$ with $\Pi=O(\sigma)$ and $\Pi(\sigma,\mu)R\equiv-R\Pi(\sigma,\mu)$. The reason is that such a term can be incorporated into $\zeta(y,z,\sigma,\mu)=O_2(y,z,\sigma)$, cf.\ equations~(2.2) in \cite{S16RCD}.
\end{rmk}

\begin{rmk}\label{Aposteriori}
Within the so-called \emph{a~posteriori} format of KAM theorems, one considers invariant tori in dynamical systems that are not assumed to be nearly integrable in any sense, see Chapter~4 of the book \cite{HCFLM16} and references therein. To the best of the author's knowledge, the a~posteriori approach to the reversible contexts has not been implemented yet.
\end{rmk}

\subsection{Aim of the Present Paper}\label{aim}

The eigenvalues of the matrix $M(\mu)$ in~\eqref{unperturbed} come in pairs $(\lambda,-\lambda)$ for each $\mu$, and generically $\det M(\mu)\neq 0$.

\begin{dfn}\label{fM}
Let a matrix $M\in\GL(2p,\mR)$ anti-commute with an involutive $2p\times 2p$ matrix with eigenvalues $1$ and $-1$ of multiplicity $p$ each. We write that the spectrum of $M$ has the form $\fM(\nu_1,\nu_2,\nu_3;\alpha,\beta)$ where $\nu_1,\nu_2,\nu_3\in\mZ_+$, $\nu_1+\nu_2+2\nu_3=p$, and $\alpha\in\mR^{\nu_1+\nu_3}$, $\beta\in\mR^{\nu_2+\nu_3}$ are two vectors with positive components, if $\det M\neq 0$ and the eigenvalues of $M$ have the form
\begin{gather*}
\pm\alpha_1,\ldots,\pm\alpha_{\nu_1}, \qquad \pm\rmi\beta_1,\ldots,\pm\rmi\beta_{\nu_2}, \\
\pm\alpha_{\nu_1+1}\pm\rmi\beta_{\nu_2+1},\ldots,\pm\alpha_{\nu_1+\nu_3}\pm\rmi\beta_{\nu_2+\nu_3}.
\end{gather*}
\end{dfn}

Assume that the spectrum of $M(\mu)$ is simple and has the form $\fM\bigl(\nu_1,\nu_2,\nu_3;\alpha(\mu),\beta(\mu)\bigr)$ for each $\mu\in\cO_s(0)$ where $\nu_1+\nu_2+2\nu_3=p$. For $\sigma=0$ and any $\mu$, the reducible invariant $n$-torus $\cT_\mu=\{y=0, \, z=0\}$ of the system~\eqref{unperturbed} has the zero Floquet exponent of multiplicity $m$ and $2p$ nonzero Floquet exponents
\begin{equation}
\begin{gathered}
\pm\alpha_j(\mu), \;\; 1\leq j\leq\nu_1; \qquad \pm\rmi\beta_j(\mu), \;\; 1\leq j\leq\nu_2; \\
\pm\alpha_{\nu_1+j}(\mu)\pm\rmi\beta_{\nu_2+j}(\mu), \;\; 1\leq j\leq\nu_3.
\end{gathered}
\label{alphabeta}
\end{equation}
For the unperturbed systems~\eqref{unperturbed}, various KAM theorems can be formulated. Let us mention four setups.

A) First, one can establish the so-called ``source'' (or Broer--Huitema--Takens-like) theorem where the frequencies $\Omega_i(\mu)$, $1\leq i\leq n$, and the nonzero Floquet exponents~\eqref{alphabeta} of the unperturbed tori are assumed to depend on $\mu$ in the ``most nondegenerate'' way, i.e., the mapping
\begin{equation}
\mu \mapsto \bigl(\Omega(\mu),\alpha(\mu),\beta(\mu)\bigr)\in\mR^{n+p}
\label{mumapsto}
\end{equation}
is submersive. This case was dealt with in our paper \cite{S16RCD}. According to the source theorem, \emph{all} the unperturbed tori $\cT_\mu$ with frequencies and nonzero Floquet exponents satisfying a suitable Diophantine condition persist under small $G$-reversible perturbations of the systems~\eqref{unperturbed}. The corresponding perturbed $n$-tori possess \emph{the same} frequency vectors and Floquet matrices and constitute a Whitney smooth family. The submersivity of the mapping~\eqref{mumapsto} is analogous to the classical Kolmogorov nondegeneracy condition for the unperturbed Lagrangian invariant tori in the Hamiltonian isotropic context without external parameters \cite{AKN06,BS10,D14}.

B) Second, one may consider weaker nondegeneracy conditions yielding just \emph{partial preservation} of the frequencies and nonzero Floquet exponents of the unperturbed tori $\cT_\mu$ under small $G$-reversible perturbations of the systems~\eqref{unperturbed}. This means that one can set up a correspondence between the unperturbed and perturbed $n$-tori in such a way that a prescribed \emph{subcollection} of the frequencies $\Omega_i(\mu)$, the positive real parts $\alpha_j(\mu)$ of the Floquet exponents~\eqref{alphabeta}, and the positive imaginary parts $\beta_j(\mu)$ of the Floquet exponents~\eqref{alphabeta} of the unperturbed tori $\cT_\mu$ coincides with the matching subcollection of the spectral characteristics of the corresponding perturbed tori. The partial preservation theorem is the subject of the present paper. In the extreme case of very weak (R\"ussmann-like \cite{R01RCD,R05RCD}) nondegeneracy conditions, perturbed systems still admit a Whitney smooth family of reducible invariant $n$-tori but it is impossible to assign the unperturbed tori to the perturbed ones in any reasonable way.

C) Third, one may examine the situation where reversible perturbations of the systems~\eqref{unperturbed} are non-autonomous and depend on time quasi-periodically with $N$ basic frequencies. In this setting studied in our paper \cite{S17MMJ}, the perturbed tori in the extended phase space are of dimension $n+N$.

D) Fourth, assuming that $\nu_2>0$, it probably makes sense to look for invariant $(n+d)$-tori $\fV^{n+d}$ ``around'' the $n$-tori $\cT_\mu$, $d=1,\ldots,\nu_2$, in the systems~\eqref{unperturbed} themselves and in their small $G$-reversible perturbations. One may speak of the excitation of the elliptic normal modes (i.e., of the purely imaginary Floquet exponents $\pm\rmi\beta_j(\mu)$, $1\leq j\leq\nu_2$) of the unperturbed tori $\cT_\mu$. This is the subject of future publications.

In the reversible context~1, Hamiltonian isotropic context, volume preserving context, and dissipative context, the four analogous setups have been more or less thoroughly explored, see the works \cite{BH95JDDE,BHS96G,BHS96LNM,BHT90,BS10,S95C,S97,S01,S06N,S07DCDS,S07TMIS} and references therein. The relevant source theorems were proven by H.W.~Broer, G.B.~Huitema, and F.~Takens in \cite{BH95JDDE,BHT90} (some generalizations are contained in the papers \cite{BCHV09PD,BHN07JDE,W10DCDSS} of Broer's group). R\"ussmann-like nondegeneracy conditions were used in \cite{BHS96G,BHS96LNM,BS10,S95C} (see also \cite{R01RCD} for R\"ussmann's original formulation in the Hamiltonian isotropic context), the general partial preservation theorems were deduced in \cite{S06N,S07TMIS}, the perturbations quasi-periodic in time were handled in \cite{S07DCDS}, and the excitation of elliptic normal modes was treated in \cite{BHS96LNM,S95C,S97,S01}.

Moreover, in all our works \cite{BHS96G,BHS96LNM,BS10,S95C,S06N,S07DCDS,S07TMIS} devoted to the four ``conventional'' KAM contexts, the results in the setups~B) and~C) were obtained as corollaries of the corresponding source theorems (whence the name ``source''). The main reduction technique we employed is called the \emph{Herman method}. This method is specifically adapted to construct invariant tori in perturbations of integrable or partially integrable systems with weak nondegeneracy conditions. It was proposed by M.R.~Herman in 1990 in his talk at an international conference on dynamical systems in Lyons (cf.\ \S~4.6.2 in \cite{Y92A}). The results in the setup~D) in \cite{BHS96LNM,S95C,S97,S01} were obtained mainly as corollaries of the results in the setup~B) with the weakest nondegeneracy conditions (so, in the long run, as corollaries of the source theorems as well). In fact, the excitation of the elliptic normal modes is only possible in the volume preserving context for tori of phase space codimension $2$ \cite{S01}, in the Hamiltonian isotropic context \cite{BHS96LNM,S97}, and in the reversible context~1 \cite{BHS96LNM,S95C}.

The idea of the Herman approach is (roughly speaking) as follows. First, by adding \emph{new} external parameters, one achieves full control of the frequencies and Floquet exponents of the unperturbed tori (the appropriate analogue of the mapping~\eqref{mumapsto} becomes submersive). The corresponding source theorem can now be applied to the new systems. Now, using the Whitney smoothness of the family of the perturbed invariant tori, the implicit function theorem, and a suitable number-theoretical lemma concerning Diophantine approximations on submanifolds of Euclidean spaces (or, as one sometimes says, Diophantine approximations of dependent quantities), one can ``extract'' information on invariant tori in the original systems (i.e., the systems without the additional external parameters). Within this procedure, all the cumbrous and laborious ``KAM machinery'' is required only to prove the source theorem and is not needed any longer as one reduces theorems with degeneracies to the source theorem.

In the reversible context~2, we also used the Herman method in the setup~C) \cite{S17MMJ} and apply this technique again in the present paper in the setup~B). So, the present paper contributes to the program \cite{S08N,S11RCD} of carrying over the results of \cite{BH95JDDE,BHS96G,BHS96LNM,BHT90,BS10,S95C,S97,S01,S06N,S07DCDS,S07TMIS} to the more involved reversible context~2 without making the proofs more complicated.

\begin{rmk}\label{original}
Partial preservation of frequencies (or frequency ratios) of the unperturbed invariant tori in the Hamiltonian isotropic context was first considered in \cite{CLY02JNS,LY05JDE,L05NA}. These papers do not employ any Herman-type reduction techniques; accordingly, the proofs in \cite{CLY02JNS,LY05JDE,L05NA} are very difficult and involve the so-called quasi-linear infinite iterative scheme.
\end{rmk}

\begin{rmk}\label{jumphop}
The codimension of the tori $\fV^{n+d}$ in the setup~D) above is equal to $m+2p-d$. Consequently, for $\frac{1}{2}(m+2p-d)\leq p=\dim\Fix G$, i.e., for $d\geq m$, the tori $\fV^{n+d}$ pertain to the reversible context~1. Thus, while examining the excitation of elliptic normal modes, one may pass from the reversible context~2 to the context~1 (cf.\ \cite{S12IM}). Similarly, while studying the destruction of unperturbed invariant tori with resonant frequencies, one may pass from the reversible context~1 to the context~2 \cite{S11RCD,S12IM,S16RCD}. Indeed, when a resonant invariant torus $\cT$ of a $G$-reversible system breaks up into a finite collection of perturbed invariant tori $\fW_1,\ldots,\fW_{\fl}$ of smaller dimension, it is possible that $\frac{1}{2}\codim\cT\leq\dim\Fix G$, but $\frac{1}{2}\codim\fW_i>\dim\Fix G$. One may suspect that the excitation of elliptic normal modes in the reversible context~2 is a much more complicated phenomenon than that in the context~1 (but it would be naive to hope that the destruction of resonant unperturbed tori is easier to study in the reversible context~2 than in the context~1). One more technique of ``moving'' from the reversible context~1 to the context~2 was developed in the paper \cite{S16RCD} where we proved the source theorem for the reversible context~2. This theorem was obtained in \cite{S16RCD} (also by Herman-like arguments) as a corollary of the main result of the article \cite{BCHV09PD} which concerns systems within the reversible context~1 with a singular normal behavior of the invariant tori.
\end{rmk}

\begin{rmk}\label{Moser}
The main tool in our first three papers \cite{S11RCD,S12MMJ,S12IM} on the reversible context~2 was the Moser modifying terms theory \cite{M67MA}.
\end{rmk}

\begin{rmk}\label{conformal}
We would like to emphasize that throughout this paper, the word ``dissipative'' means ``relating to no structure on the phase space''. For instance, conformally Hamiltonian vector fields $V$ and conformally symplectic diffeomorphisms $A$ which have been extensively studied lately in KAM theory (see \cite{CCdL17N} and references therein) are defined by the identities $d(i_V\bfomega)\equiv\eta\bfomega$ and $A^\ast\bfomega\equiv\pm\rme^\eta\bfomega$ with constant $\eta\neq 0$ and are therefore \emph{not} dissipative in this sense ($\bfomega$ being a symplectic structure on the phase space). However, conformally Hamiltonian systems are dissipative in another sense, namely, their dynamics exhibits no conservative patterns.
\end{rmk}

Like in our previous works on the ``conventional'' KAM contexts and the reversible context~2, we confine ourselves with analytic systems but there is no doubt that our results (Theorems~\ref{thmain} and~\ref{Ruessmann} below) can be carried over to Gevrey regular or merely infinitely differentiable systems and even to $C^r$-smooth systems for finite (but sufficiently large) $r$. Similarly, the families of analytic perturbed invariant tori in Theorems~\ref{thmain}, \ref{Ruessmann}, and~\ref{thsource} below are claimed to be $C^\infty$-smooth in the sense of Whitney, but these families are certainly Gevrey regular in the sense of Whitney (cf.\ \cite{W10DCDSS}).

The paper is organized as follows. In Section~\ref{numbertheory}, we formulate the Diophantine lemma (Lemma~\ref{lmDiophantine}) to be used in the Herman procedure. The main result of the paper (Theorem~\ref{thmain}) is stated in Section~\ref{mainresult}. In Section~\ref{source}, we give a precise formulation of the source theorem for the reversible context~2 (Theorem~\ref{thsource}) in the form we need. A proof of the main result is presented in Section~\ref{proof}. Finally, in Section~\ref{multiparameter}, we give a rigorous proof of the fact that the occurrence of invariant tori in the reversible context~2 requires many external parameters.

\section{Diophantine Lemma}\label{numbertheory}

\begin{dfn}[\cite{S06N,S07DCDS,S07TMIS,S17MMJ}]\label{Diophantine}
Let $n\in\mZ_+$ and $\nu\in\mZ_+$. Given $\tau\geq 0$, $\gamma>0$, and $L\in\mN$, a pair of vectors
\begin{equation}
\Omega\in\mR^n, \quad \beta\in\mR^\nu
\label{nnu}
\end{equation}
is said to be \emph{affinely $(\tau,\gamma,L)$-Diophantine}, if the inequality
\[
\bigl|\langle\Omega,k\rangle + \langle\beta,\ell\rangle\bigr| \geq \gamma|k|^{-\tau}
\]
holds for any $k\in\mZ^n\setminus\{0\}$ and $\ell\in\mZ^\nu$ such that $|\ell|\leq L$.
\end{dfn}

Clearly, if $n\in\mN$ and a pair of vectors~\eqref{nnu} is affinely $(\tau,\gamma,L)$-Diophantine, then the vector $\Omega\in\mR^n$ is $(\tau,\gamma)$-Diophantine in the usual sense, so that $\tau\geq n-1$. If $n=0$ then a pair of vectors~\eqref{nnu} is affinely $(\tau,\gamma,L)$-Diophantine for any $\tau$, $\gamma$, $L$, $\nu$, and $\beta\in\mR^\nu$ \cite{S07TMIS}.

\begin{dfn}[\cite{S06N,S07DCDS,S07TMIS,S17MMJ}]\label{nondegenerate}
Let $s\in\mN$, $n\in\mZ_+$, and $\nu\in\mZ_+$. We will adopt the standard multi-index notation
\[
q!=q_1!q_2!\cdots q_s!, \quad
\mu^q=\mu_1^{q_1}\mu_2^{q_2}\cdots\mu_s^{q_s}, \quad
D_\mu^q\Omega=\frac{\partial^{|q|}\Omega}{\partial\mu_1^{q_1}\partial\mu_2^{q_2}\cdots\partial\mu_s^{q_s}},
\]
where $q\in\mZ_+^s$, $\mu\in\mR^s$, and $\Omega$ is a (vector-valued) function $C^{|q|}$-smooth in $\mu$. Let $\fA\subset\mR^s$ be an open domain and let $Q\in\mN$, $L\in\mN$. Consider a pair of $C^Q$-smooth mappings $\Omega:\fA\to\mR^n$, $\beta:\fA\to\mR^\nu$. If $n>0$, introduce the notation
\[
\rho^Q(\mu) = \min_{\|e\|=1} \max_{J=1}^Q J! \max_{\|\fu\|=1}
\left|\sum_{|q|=J} \bigl\langle D_\mu^q\Omega(\mu),e\bigr\rangle \frac{\fu^q}{q!}\right|
\]
($q\in\mZ_+^s$, $e\in\mR^n$, $\fu\in\mR^s$) for $\mu\in\fA$. If $\nu>0$, introduce the notation
\[
\kappa_\ell^Q(\mu) = \max_{J=1}^Q J! \max_{\|\fu\|=1}
\left|\sum_{|q|=J} \bigl\langle D_\mu^q\beta(\mu),\ell\bigr\rangle \frac{\fu^q}{q!}\right|
\]
($q\in\mZ_+^s$, $\fu\in\mR^s$) for $\mu\in\fA$, $\ell\in\mZ^\nu$. The pair of mappings $\Omega$, $\beta$ is said to be \emph{affinely $(Q,L)$-nondegenerate} at a point $\mu\in\fA$ if one of the following four conditions is satisfied.

1) $n>0$, $\nu>0$, $\rho^Q(\mu)>0$, and
\[
\max_{1\leq |q|\leq Q} \Bigl| \bigl\langle D_\mu^q\Omega(\mu),k\bigr\rangle + \bigl\langle D_\mu^q\beta(\mu),\ell\bigr\rangle \Bigr| > 0
\]
($q\in\mZ_+^s$) for all $k\in\mZ^n$ and $\ell\in\mZ^\nu$ such that $1\leq |\ell|\leq L$ and $\|k\| \leq \kappa_\ell^Q(\mu)\big/\rho^Q(\mu)$.

2) $n>0$, $\nu=0$, and $\rho^Q(\mu)>0$.

3) $n=0$, $\nu>0$, and $\kappa_\ell^Q(\mu)>0$ for all $\ell\in\mZ^\nu$ such that $1\leq |\ell|\leq L$.

4) $n=\nu=0$.
\end{dfn}

Note that for any (vector-valued) $C^J$-smooth function $H$ defined in $\fA\subset\mR^s$ ($J\in\mZ_+$) and any $\mu\in\fA$, $\fu\in\mR^s$ one has
\[
J! \sum_{|q|=J} D_\mu^qH(\mu)\frac{\fu^q}{q!} =
\left. \frac{d^J}{dt^J}H(\mu+t\fu) \right|_{t=0}
\]
($q\in\mZ_+^s$). The inequality $\rho^Q(\mu)>0$ (for $n>0$) means that the collection of all the $\binom{s+Q}{s}-1$ partial derivatives of $\Omega$ at $\mu$ of all the orders from $1$ to $Q$ spans $\mR^n$, i.e., the linear hull of these derivatives is $\mR^n$ (a R\"ussmann-type property \cite{R01RCD}). The inequality $\kappa_\ell^Q(\mu)>0$ (for $\nu>0$ and some $\ell\in\mZ^\nu\setminus\{0\}$) means that at least one of the $\binom{s+Q}{s}-1$ partial derivatives of $\beta$ at $\mu$ of all the orders from $1$ to $Q$ is not orthogonal to $\ell$. It is not hard to verify that if a pair of mappings $\Omega$, $\beta$ is affinely $(Q,L)$-nondegenerate at a point $\mu\in\fA$, then it is affinely $(Q,L)$-nondegenerate at any point $\mu'\in\fA$ sufficiently close to $\mu$.

\begin{lmm}[\cite{S07TMIS}]\label{lmDiophantine}
Let $s\in\mN$, $n\in\mZ_+$, $\nu\in\mZ_+$, $Q\in\mN$, and $L\in\mN$. Let also $\fA\subset\mR^s$ be an open domain, $A\subset\fA$ a subset of $\fA$ diffeomorphic to a closed $s$-dimensional ball, and $B$ an arbitrary compact metric space. Suppose that mappings
\[
\Omega:\fA\times B\to\mR^n \quad \text{and} \quad \beta:\fA\times B\to\mR^\nu
\]
are $C^Q$-smooth in $a\in\fA$ and, moreover, all the partial derivatives of the functions $\Omega$ and $\beta$ with respect to $a_1,\ldots,a_s$ of any order from $1$ to $Q$ are continuous in $(a,b)\in\fA\times B$ \textup{(}rather than only in $a\in\fA$\textup{)}. Let the pair of mappings
\begin{equation}
a\mapsto\Omega(a,b)\in\mR^n \quad \text{and} \quad a\mapsto\beta(a,b)\in\mR^\nu
\label{abstractpair}
\end{equation}
be affinely $(Q,L)$-nondegenerate at each point $a\in A$ for any fixed value of $b\in B$. Then

\textup{(1)} there exists a number $\delta>0$ and

\textup{(2)} for every $n\add\in\mZ_+$, $\nu\add\in\mZ_+$, $\tau_\ast\geq\max(0,n\add-1)$, $\gamma_\ast>0$, $\vare\in(0,1)$ and every $\tau$ such that $\tau>(n+n\add)Q$ and $\tau\geq\tau_\ast$, there exists a number $\gamma=\gamma_0(\vare,\tau,\gamma_\ast)>0$ such that the following holds. Let
\[
\tOmega:\fA\times B\to\mR^n \quad \text{and} \quad \tbeta:\fA\times B\to\mR^\nu
\]
be any mappings $C^Q$-smooth in $a\in\fA$ and such that all the partial derivatives of each component of the differences $\tOmega-\Omega$ and $\tbeta-\beta$ with respect to $a_1,\ldots,a_s$ of any order from $1$ to $Q$ are smaller than $\delta$ in absolute value everywhere in $\fA\times B$. Let
\begin{equation}
\Omega\add:B\to\mR^{n\add} \quad \text{and} \quad \beta\add:B\to\mR^{\nu\add}
\label{additionalpair}
\end{equation}
be arbitrary mappings. Then, for any $b\in B$ such that the pair of vectors $\Omega\add(b)$, $\beta\add(b)$ is affinely $(\tau_\ast,\gamma_\ast,L)$-Diophantine, the Lebesgue measure of the set of those points $a\in A$ for which the pair of vectors
\[
\Bigl( \tOmega(a,b),\Omega\add(b) \Bigr)\in\mR^{n+n\add}, \quad \Bigl( \tbeta(a,b),\beta\add(b) \Bigr)\in\mR^{\nu+\nu\add}
\]
is affinely $(\tau,\gamma,L)$-Diophantine, is greater than $(1-\vare)\meas_sA$.
\end{lmm}

Here and henceforth, $\meas_s$ denotes the Lebesgue measure in $\mR^s$. Some particular cases of Lemma~\ref{lmDiophantine} are formulated in \cite{S06N,S07DCDS,S17MMJ}.

\begin{xpl}\label{compact}
The compactness of $B$ in Lemma~\ref{lmDiophantine} is essential. For instance, suppose that $n\in\mN$ and a pair of $C^Q$-smooth mappings
\begin{equation}
\Omega_0:\fA\to\mR^n \quad \text{and} \quad \beta_0:\fA\to\mR^\nu
\label{reducedpair}
\end{equation}
is affinely $(Q,L)$-nondegenerate at each point $a\in A$. Let $B=[1,+\infty)$ and $\Omega(a,b)=\Omega_0(a)/b$, $\beta(a,b)=\beta_0(a)/b$. The pair of mappings~\eqref{abstractpair} is affinely $(Q,L)$-nondegenerate at each point $a\in A$ for any fixed value of $b\in B$. Assume also that all the partial derivatives of each component of the functions~\eqref{reducedpair} of any order from $1$ to $Q$ are no greater than a certain number $\fD<+\infty$ in absolute value everywhere in $\fA$. Given $\delta>0$, let $c_1=\max(\fD/\delta,1)$ and choose an arbitrary number $c_2>c_1$. Consider an arbitrary function $\vart:B\to\mR$ such that $\vart(b)=1$ for $1\leq b\leq c_1$, $0<\vart(b)<1$ for $c_1<b<c_2$, and $\vart(b)=0$ for $b\geq c_2$ (such a function can be chosen to be $C^\infty$-smooth, but we will not use this fact). Set
\begin{equation}
\tOmega(a,b)=\vart(b)\Omega(a,b)=\vart(b)\Omega_0(a)/b \quad \text{and} \quad \tbeta(a,b)=\vart(b)\beta(a,b)=\vart(b)\beta_0(a)/b.
\label{perturbedpair}
\end{equation}
Since
\[
\frac{\fD\bigl(1-\vart(b)\bigr)}{b} < \delta
\]
for any $b\in B$, all the partial derivatives of each component of the differences $\tOmega-\Omega$ and $\tbeta-\beta$ with respect to $a$ of any order from $1$ to $Q$ are smaller than $\delta$ in absolute value everywhere in $\fA\times B$. Now let $n\add=\nu\add=0$, i.e., let the mappings~\eqref{additionalpair} be absent. Given arbitrary $\vare\in(0,1)$, $\tau>nQ$, $\gamma>0$, one \emph{cannot} assert that for any $b\in B$ the Lebesgue measure of the set $\cA_b$ of those points $a\in A$ for which the pair of vectors~\eqref{perturbedpair} is affinely $(\tau,\gamma,L)$-Diophantine, is greater than $(1-\vare)\meas_sA$. Indeed, $\cA_b$ is empty for $b\geq c_2$ because both the vectors~\eqref{perturbedpair} are zero for each $a$ if $b\geq c_2$.
\end{xpl}

\section{The Main Result}\label{mainresult}

In Sections~\ref{mainresult} and~\ref{source}, we will sometimes write $\{0\in\mR^s\}$ instead of $\{0\}$ with $0\in\mR^s$.

Let $n\in\mZ_+$, $m\in\mN$, $p\in\mZ_+$, $s\in\mN$. Consider an analytic $(m+s)$-parameter family of analytic differential equations
\begin{equation}
\begin{aligned}
\dot{x} &= \Omega(\mu)+\Delta(\sigma,\mu)+\xi(y,z,\sigma,\mu)+f(x,y,z,\sigma,\mu), \\
\dot{y} &= \sigma+\eta(y,z,\sigma,\mu)+g(x,y,z,\sigma,\mu), \\
\dot{z} &= M(\mu)z+\zeta(y,z,\sigma,\mu)+h(x,y,z,\sigma,\mu),
\end{aligned}
\label{perturbed}
\end{equation}
where $x\in\mT^n$, $y\in\cO_m(0)$, $z\in\cO_{2p}(0)$ are the phase space variables, $\sigma\in\cO_m(0)$ and $\mu\in\cO_s(0)$ are external parameters, $M$ is a $2p\times 2p$ matrix-valued function, $\Delta=O(\sigma)$ and $\xi=O(y,z)$, $\eta=O_2(y,z)$, $\zeta=O_2(y,z,\sigma)$. The functions $\Omega$, $M$, $\Delta$, $\xi$, $\eta$, $\zeta$ are supposed to be fixed whereas the terms $f$, $g$, $h$ are small perturbations, cf.~\eqref{unperturbed}. Let the systems~\eqref{perturbed} be reversible with respect to the phase space involution~\eqref{involution}, where $R\in\GL(2p,\mR)$ is an involutive matrix with eigenvalues $1$ and $-1$ of multiplicity $p$ each, $M(\mu)R\equiv-RM(\mu)$, and the spectrum of $M(0)$ is simple. One may assume that the spectrum of $M(\mu)$ is simple for each $\mu$ and has the form $\fM\bigl(\nu_1,\nu_2,\nu_3;\alpha(\mu),\beta(\mu)\bigr)$ where $\nu_1+\nu_2+2\nu_3=p$ (see Definition~\ref{fM}). Introduce the notation $\nu=\nu_2+\nu_3\in\mZ_+$.

Choose arbitrary (possibly, empty) subsets of indices
\begin{gather*}
\fS_1\subset\{1;2;\ldots;n\}, \quad \fS_2\subset\{1;2;\ldots;\nu_1+\nu_3\}, \quad \fS_3\subset\{1;2;\ldots;\nu\}, \\
\fT\subset\{1;2;\ldots;s\}
\end{gather*}
such that
\[
0\leq\#\fS_1+\#\fS_2+\#\fS_3=\#\fT\leq\min(n+p,s-1).
\]
Here and henceforth, $\#$ denotes the number of elements of a finite set. We are interested in the preservation of the frequencies $\Omega_i$ (of the unperturbed invariant tori $\cT_\mu=\{y=0, \, z=0\}$ at $\sigma=0$) with $i\in\fS_1$, the real parts $\alpha_j$ of the Floquet exponents with $j\in\fS_2$, and the imaginary parts $\beta_j$ of the Floquet exponents with $j\in\fS_3$. We will write
\begin{gather}
\Omega_+=(\Omega_i \mid i\in\fS_1), \quad \Omega_-=(\Omega_i \mid i\notin\fS_1), \nonumber \\
\alpha_+=(\alpha_j \mid j\in\fS_2), \quad \alpha_-=(\alpha_j \mid j\notin\fS_2), \nonumber \\
\beta_+=(\beta_j \mid j\in\fS_3), \quad \beta_-=(\beta_j \mid j\notin\fS_3), \nonumber \\
\mu_+=(\mu_l \mid l\in\fT), \quad \mu_-=(\mu_l \mid l\notin\fT);
\label{fT}
\end{gather}
similar notation will be used below without special mention for vector quantities denoted by the letters $\Omega$, $\alpha$, $\beta$, $\mu$ with superscripts or diacritical marks. Set
\[
\#\fS_1=d_1, \quad \#\fS_2=d_2, \quad \#\fS_3=d_3, \qquad d_1+d_2+d_3=d=\#\fT,
\]
so that
\begin{gather}
0\leq d_1\leq n, \quad 0\leq d_2\leq\nu_1+\nu_3, \quad 0\leq d_3\leq\nu, \nonumber \\
0\leq d\leq\min(n+p,s-1).
\label{rangeofd}
\end{gather}
For any vector $b\in\mR^d$, we will write
\[
b^{:1}=(b_1,\ldots,b_{d_1})\in\mR^{d_1}, \quad
b^{:2}=(b_{d_1+1},\ldots,b_{d_1+d_2})\in\mR^{d_2}, \quad
b^{:3}=(b_{d_1+d_2+1},\ldots,b_d)\in\mR^{d_3}.
\]
We will also use the notation
\[
\fP_0=\bigl(\Omega_+(0),\alpha_+(0),\beta_+(0)\bigr)\in\mR^d.
\]

\begin{thm}\label{thmain}
Suppose that either

$d=0$

\noindent or

$d>0$ and the Jacobian
\begin{equation}
\frac{\partial(\Omega_+,\alpha_+,\beta_+)}{\partial\mu_+}
\label{Jacobian}
\end{equation}
of order $d$ does not vanish at $\mu=0$. This implies, in particular, that $(\Omega_+,\alpha_+,\beta_+)$ can be used as a part of a new coordinate frame near the origin of $\mR^s$. In other words, there exists an analytic change of coordinates $\mu=\mu(a,b)$ in a neighborhood of $\mu=0$ such that
\[
a\in\cO_{s-d}(0), \quad b\in\cO_d(\fP_0), \qquad \mu(0,\fP_0)=0,
\]
and
\[
(\Omega_+,\alpha_+,\beta_+)\big|_{\mu=\mu(a,b)}\equiv b,
\]
more precisely,
\begin{equation}
\Omega_+\bigl(\mu(a,b)\bigr)\equiv b^{:1}, \quad
\alpha_+\bigl(\mu(a,b)\bigr)\equiv b^{:2}, \quad
\beta_+\bigl(\mu(a,b)\bigr)\equiv b^{:3}.
\label{rasklad}
\end{equation}
Assume also that a change of coordinates $\mu=\mu(a,b)$ with this property can be chosen in such a way that the pair of mappings
\begin{equation}
a\mapsto\Omega_-\bigl(\mu(a,0)\bigr)\in\mR^{n-d_1} \quad \text{and} \quad a\mapsto\beta_-\bigl(\mu(a,0)\bigr)\in\mR^{\nu-d_3}
\label{pair}
\end{equation}
is affinely $(Q,2)$-nondegenerate at $a=0$ for some number $Q\in\mN$ \textup{(}see Definition~\ref{nondegenerate}\textup{)}.

Then there exist a closed $(s-d)$-dimensional ball $A\subset\mR^{s-d}$ centered at the origin and a closed $d$-dimensional ball $B\subset\mR^d$ centered at the point $\fP_0$ such that the following holds. Set
\begin{equation}
\Gamma=\bigl\{ \mu(a,b) \bigm| a\in A, \, b\in B \bigr\}\subset\mR^s
\label{Gamma}
\end{equation}
\textup{(}$0\in\Gamma$\textup{)}. Then for every complex neighborhood
\begin{equation}
\cC\subset(\mC/2\pi\mZ)^n\times\mC^{2m+2p+s}
\label{neighborhood}
\end{equation}
of the set
\begin{equation}
\mT^n\times\{0\in\mR^m\}\times\{0\in\mR^{2p}\}\times\{0\in\mR^m\}\times\Gamma,
\label{ostov}
\end{equation}
every $\cL\in\mN$, $\vare_1>0$, $\vare_2\in(0,1)$, $\vare_3\in(0,1)$, $\tau_\ast\geq\max(0,d_1-1)$, $\gamma_\ast>0$, and every $\tau$ such that $\tau>nQ$ and $\tau\geq\tau_\ast$, there are numbers $\delta>0$ and $\gamma\in(0,\gamma_\ast]$ with the following properties.

Suppose that the perturbation terms $f$, $g$, $h$ in~\eqref{perturbed} can be holomorphically continued to the neighborhood $\cC$ and $|f|<\delta$, $|g|<\delta$, $|h|<\delta$ in $\cC$. Consider the closed $(s-d)$-dimensional ball $\tA\subset A$ centered at the origin and the closed $d$-dimensional ball $\tB\subset B$ centered at the point $\fP_0$ such that
\begin{equation}
\meas_{s-d}\tA = (1-\vare_3)\meas_{s-d}A, \quad \meas_d\tB = (1-\vare_3)\meas_dB
\label{tAtB}
\end{equation}
and set
\begin{equation}
\tGamma=\Bigl\{ \mu(a,b) \Bigm| a\in\tA, \, b\in\tB \Bigr\}\subset\Gamma
\label{tGamma}
\end{equation}
\textup{(}$0\in\tGamma$\textup{)}. Then there exist functions
\begin{equation}
\begin{gathered}
\Theta:\tGamma\to\mR^m, \quad \Xi:\tGamma\to\mR^s, \\
\tOmega:\tGamma\to\mR^n, \quad \tM:\tGamma\to\gl(2p,\mR)
\end{gathered}
\label{functions}
\end{equation}
and a change of variables
\begin{equation}
\begin{aligned}
x &= \bx+X(\bx,\mu), \\
y &= \by+Y^0(\bx,\mu)+Y^1(\bx,\mu)\by+Y^2(\bx,\mu)\bz, \\
z &= \bz+Z^0(\bx,\mu)+Z^1(\bx,\mu)\by+Z^2(\bx,\mu)\bz
\end{aligned}
\label{change}
\end{equation}
for each $\mu\in\tGamma$ with $\bx\in\mT^n$, $\by\in\cO_m(0)$, $\bz\in\cO_{2p}(0)$ such that the following is valid.

\textup{i)} The functions~\eqref{functions} are $C^\infty$-smooth, and all the partial derivatives of each component of the functions $\Theta$, $\Xi$, $\tOmega-\Omega$, $\tM-M$ of any order from $0$ to $\cL$ are smaller than $\vare_1$ in absolute value everywhere in $\tGamma$. If $d=0$ then $\Xi\equiv 0$. The coefficients $X$, $Y^0$, $Y^1$, $Y^2$, $Z^0$, $Z^1$, $Z^2$ in~\eqref{change} are mappings ranging in $\mR^n$, $\mR^m$, $\gl(m,\mR)$, $\mR^{m\times 2p}$, $\mR^{2p}$, $\mR^{2p\times m}$, $\gl(2p,\mR)$, respectively. These mappings are analytic in $\bx$ and $C^\infty$-smooth in $\mu$. All the partial derivatives of each component of these mappings of any order from $0$ to $\cL$ are smaller than $\vare_1$ in absolute value everywhere in $\mT^n\times\tGamma$.

\textup{ii)} For each $\mu\in\tGamma$, the change of variables~\eqref{change} commutes with the involution~\eqref{involution} in the sense that in the new variables $(\bx,\by,\bz)$, the involution $G$ takes the form
\[
G:(\bx,\by,\bz)\mapsto(-\bx,-\by,R\bz).
\]
There holds the identity $\tM(\mu)R\equiv-R\tM(\mu)$.

\textup{iii)} The spectrum of $\tM(\mu)$ is simple and has the form $\fM\Bigl(\nu_1,\nu_2,\nu_3;\talpha(\mu),\tbeta(\mu)\Bigr)$ for each $\mu\in\tGamma$ \textup{(}see Definition~\ref{fM}\textup{)}, and the identities
\begin{equation}
\tOmega_+\equiv\Omega_+, \quad \talpha_+\equiv\alpha_+, \quad \tbeta_+\equiv\beta_+
\label{identities}
\end{equation}
are valid in $\tGamma$.

\textup{iv)} For any point $b\in\tB$ such that the pair of vectors $b^{:1}\in\mR^{d_1}$, $b^{:3}\in\mR^{d_3}$ is affinely $(\tau_\ast,\gamma_\ast,2)$-Diophantine \textup{(}see Definition~\ref{Diophantine}\textup{)}, there exists a set $\cG_b\subset\tA$ that satisfies the following conditions.

\textup{(a)} $\meas_{s-d}\cG_b > (1-\vare_2)\meas_{s-d}\tA$.

\textup{(b)} For any point $a\in\cG_b$, the pair of vectors $\tOmega(\mu^0)\in\mR^n$, $\tbeta(\mu^0)\in\mR^\nu$ is affinely $(\tau,\gamma,2)$-Diophantine, where $\mu^0=\mu(a,b)$.

\textup{(c)} For any point $a\in\cG_b$, the perturbed system~\eqref{perturbed} with $\mu=\mu^0+\Xi(\mu^0)$ and $\sigma=\Theta(\mu^0)$ takes the form
\begin{equation}
\dot{\bx}=\tOmega(\mu^0)+O(\by,\bz), \quad \dot{\by}=O_2(\by,\bz), \quad \dot{\bz}=\tM(\mu^0)\bz+O_2(\by,\bz)
\label{goal}
\end{equation}
after the coordinate change~\eqref{change} with $\mu=\mu^0$.
\end{thm}

\begin{rmk}\label{tB}
One can easily show that for $d_1\in\mN$ the Lebesgue measure $\meas_d$ of the set of those points $b\in\tB$ for which the pair of vectors $b^{:1}\in\mR^{d_1}$, $b^{:3}\in\mR^{d_3}$ is \emph{not} affinely $(\tau_\ast,\gamma_\ast,2)$-Diophantine, tends to $0$ as $\gamma_\ast\to 0$ for any fixed $\tau_\ast>d_1-1$. For $d_1=0$ this set is empty for any $\tau_\ast\geq 0$ and $\gamma_\ast>0$.
\end{rmk}

So, consider an arbitrary point $\mu^\star=\mu(a^\star,b^\star)\in\tGamma$ such that the pair of vectors $b^{\star:1}=\Omega_+(\mu^\star)\in\mR^{d_1}$, $b^{\star:3}=\beta_+(\mu^\star)\in\mR^{d_3}$ is affinely $(\tau_\ast,\gamma_\ast,2)$-Diophantine, see~\eqref{rasklad}. The unperturbed invariant $n$-tori $\cT_\mu$, $\mu\in\Gamma$, such that
\begin{equation}
\Omega_+(\mu)=\Omega_+(\mu^\star), \quad \alpha_+(\mu)=\alpha_+(\mu^\star), \quad \beta_+(\mu)=\beta_+(\mu^\star)
\label{thesame}
\end{equation}
constitute an $(s-d)$-parameter smooth family: the equalities~\eqref{thesame} are equivalent to that $\mu=\mu(a,b^\star)$ for some $a\in A$. Now choose any $a\in\cG_{b^\star}$ and denote $\mu(a,b^\star)$ by $\mu^0$. The perturbed system~\eqref{perturbed} with the \emph{shifted} parameter values $\mu=\mu^0+\Xi(\mu^0)$, $\sigma=\Theta(\mu^0)$ and the involution~\eqref{involution} admit a common reducible invariant $n$-torus $\{\by=0, \, \bz=0\}$ with frequency vector $\tOmega(\mu^0)$ and Floquet matrix $\bfzero_m\oplus\tM(\mu^0)$, see~\eqref{goal}. For the frequencies $\tOmega_i(\mu^0)$ of this torus, the positive real parts $\talpha_j(\mu^0)$ of the Floquet exponents, and the positive imaginary parts $\tbeta_j(\mu^0)$ of the Floquet exponents, one has
\begin{gather*}
\tOmega_+(\mu^0)=\Omega_+(\mu^0)=\Omega_+(\mu^\star), \\
\talpha_+(\mu^0)=\alpha_+(\mu^0)=\alpha_+(\mu^\star), \\
\tbeta_+(\mu^0)=\beta_+(\mu^0)=\beta_+(\mu^\star)
\end{gather*}
according to~\eqref{identities} and~\eqref{thesame}. All these perturbed tori constitute an $(s-d)$-parameter Cantor family (the parameter being $a\in\cG_{b^\star}$). The torus $\{\by=0, \, \bz=0\}$ is analytic and depends on $(a,b^\star)$ in a $C^\infty$-way in the sense of Whitney.

Such partial preservation of the frequencies and the real and imaginary parts of the Floquet exponents of the unperturbed tori $\cT_\mu$ is essentially provided by two nondegeneracy conditions: a Broer--Huitema--Takens-type condition \cite{BH95JDDE,BHT90} on the components $\Omega_+$, $\alpha_+$, $\beta_+$ to be preserved (the Jacobian~\eqref{Jacobian} does not vanish for $\mu\in\mR^s$ near $0$) and a R\"ussmann-type condition \cite{R01RCD} on the components $\Omega_-$, $\beta_-$ (the pair of mappings~\eqref{pair} is affinely $(Q,2)$-nondegenerate for $a\in\mR^{s-d}$ near $0$). The second condition requires $s-d$ to be positive. This is the reason why we assume that $d\leq s-1$ in Theorem~\ref{thmain}, although the bound $d\leq\min(n+p,s)$ may seem more ``natural'' in~\eqref{rangeofd} than $d\leq\min(n+p,s-1)$.

In the coordinate transformation~\eqref{change}, the terms $X$, $Y^0$, and $Z^0$ are responsible for the invariance of the torus $\{\by=0, \, \bz=0\}$ while the terms $Y^1\by$, $Y^2\bz$, $Z^1\by$, and $Z^2\bz$ are responsible for its reducibility, i.e., for the variational equation along $\{\by=0, \, \bz=0\}$, see a detailed discussion in \cite{S17MMJ}.

\begin{rmk}\label{heuristic}
The frequency vectors of the unperturbed invariant tori $\cT_\mu=\{y=0, \, z=0\}$ of the systems~\eqref{unperturbed} are $\Omega(\mu)$, and generically some of these tori are resonant, some are not. Theorem~\ref{thmain} (even with $d=0$) shows that a small generic $G$-reversible perturbation of these systems preserves the family of tori $\cT_\mu$ but makes it Cantor-like (for $n\geq 2$ under mild nondegeneracy conditions). Thus, the reversible context~2 obeys the heuristic principle formulated in Section~2 in the paper \cite{BHS96G} and in \S~1.4.1 in the book \cite{BHS96LNM}.
\end{rmk}

In the four ``conventional'' KAM contexts (the reversible context~1, Hamiltonian isotropic context, volume preserving context, and dissipative context), we had the following picture \cite{BHS96G,BHS96LNM,S95C,S06N,S07TMIS}: if the differential equations depend on $\fc$ external parameters and the reducible invariant tori constitute an $\fs$-parameter Cantor-like family in the product of the phase space and the space of external parameters, then always $\fs\geq\fc$ and each torus has $\fs-\fc$ zero Floquet exponents (if $\fs=\fc$ then $\fc$ is required to be at least $1$). In the reversible context~2, on the contrary, there holds the inequality $\fc>\fs$, and each torus has $\fc-\fs$ zero Floquet exponents. Indeed, in the framework of Theorem~\ref{thmain}, $\fc=m+s$, $\fs=s$, and each torus possesses a zero Floquet exponent of multiplicity $m$. In all the five contexts, each perturbed torus has $|\fs-\fc|$ zero Floquet exponents.

\begin{rmk}\label{Brouwer}
In the Hamiltonian isotropic context and the reversible context~1, there are results on the frequency preservation where the nondegeneracy condition is formulated in terms of the Brouwer topological degree rather than in terms of the rank of a certain Jacobi-type matrix. The relevant references for reversible systems are \cite{KX14AMC,WXZ10DCDSB,WXZ16AAM,WXZ17DCDS}. As to Hamiltonian systems, we confine ourselves with the papers \cite{XL16RCD,ZXW16ANS} (see also references therein). Of all these works, the papers \cite{KX14AMC,WXZ10DCDSB,WXZ17DCDS,ZXW16ANS} employ the Herman approach. Some sets of nondegeneracy conditions in the Hamiltonian isotropic context and the reversible context~1 are reviewed in \cite{H11IM}.
\end{rmk}

For $d=0$ none of the frequencies and Floquet exponents of the unperturbed tori in Theorem~\ref{thmain} is required to be preserved (a R\"ussmann-like situation \cite{R01RCD,R05RCD}). The simplest case where $d=0$ and $p=0$ was examined in Section~5 in the paper \cite{S16RCD}. Since the case of zero $d$ is very important, we present it as a separate theorem. Consider again the system~\eqref{perturbed} of differential equations.

\begin{thm}\label{Ruessmann}
Let the pair of mappings
\[
\mu\mapsto\Omega(\mu)\in\mR^n \quad \text{and} \quad \mu\mapsto\beta(\mu)\in\mR^\nu
\]
be affinely $(Q,2)$-nondegenerate at $\mu=0$ for some number $Q\in\mN$ \textup{(}see Definition~\ref{nondegenerate}\textup{)}. Then there exists a closed $s$-dimensional ball $\Gamma\subset\mR^s$ centered at the origin and such that the following holds. For every complex neighborhood~\eqref{neighborhood} of the set~\eqref{ostov} and every $\cL\in\mN$, $\vare_1>0$, $\vare_2\in(0,1)$, $\vare_3\in(0,1)$, $\tau>nQ$, there are numbers $\delta>0$ and $\gamma>0$ with the following properties.

Suppose that the perturbation terms $f$, $g$, $h$ in~\eqref{perturbed} can be holomorphically continued to the neighborhood $\cC$ and $|f|<\delta$, $|g|<\delta$, $|h|<\delta$ in $\cC$. Consider the closed $s$-dimensional ball $\tGamma\subset\Gamma$ centered at the origin and such that
\[
\meas_s\tGamma = (1-\vare_3)\meas_s\Gamma.
\]
Then there exist functions
\begin{equation}
\Theta:\tGamma\to\mR^m, \quad \tOmega:\tGamma\to\mR^n, \quad \tM:\tGamma\to\gl(2p,\mR)
\label{withoutXi}
\end{equation}
and a change of variables~\eqref{change} for each $\mu\in\tGamma$ with $\bx\in\mT^n$, $\by\in\cO_m(0)$, $\bz\in\cO_{2p}(0)$ such that the following is valid.

\textup{i)} The functions~\eqref{withoutXi} are $C^\infty$-smooth, and all the partial derivatives of each component of the functions $\Theta$, $\tOmega-\Omega$, $\tM-M$ of any order from $0$ to $\cL$ are smaller than $\vare_1$ in absolute value everywhere in $\tGamma$. The coefficients $X$, $Y^0$, $Y^1$, $Y^2$, $Z^0$, $Z^1$, $Z^2$ in~\eqref{change} are mappings ranging in $\mR^n$, $\mR^m$, $\gl(m,\mR)$, $\mR^{m\times 2p}$, $\mR^{2p}$, $\mR^{2p\times m}$, $\gl(2p,\mR)$, respectively. These mappings are analytic in $\bx$ and $C^\infty$-smooth in $\mu$. All the partial derivatives of each component of these mappings of any order from $0$ to $\cL$ are smaller than $\vare_1$ in absolute value everywhere in $\mT^n\times\tGamma$.

\textup{ii)} For each $\mu\in\tGamma$, the change of variables~\eqref{change} commutes with the involution~\eqref{involution}. There holds the identity $\tM(\mu)R\equiv-R\tM(\mu)$.

\textup{iii)} The spectrum of $\tM(\mu)$ is simple and has the form $\fM\Bigl(\nu_1,\nu_2,\nu_3;\talpha(\mu),\tbeta(\mu)\Bigr)$ for each $\mu\in\tGamma$ \textup{(}see Definition~\ref{fM}\textup{)}.

\textup{iv)} There exists a set $\cG\subset\tGamma$ that satisfies the following conditions.

\textup{(a)} $\meas_s\cG > (1-\vare_2)\meas_s\tGamma$.

\textup{(b)} For any point $\mu\in\cG$, the pair of vectors $\tOmega(\mu)\in\mR^n$, $\tbeta(\mu)\in\mR^\nu$ is affinely $(\tau,\gamma,2)$-Diophantine \textup{(}see Definition~\ref{Diophantine}\textup{)}.

\textup{(c)} For any point $\mu\in\cG$, the perturbed system~\eqref{perturbed} with $\sigma=\Theta(\mu)$ takes the form
\begin{equation}
\dot{\bx}=\tOmega(\mu)+O(\by,\bz), \quad \dot{\by}=O_2(\by,\bz), \quad \dot{\bz}=\tM(\mu)\bz+O_2(\by,\bz)
\label{withoutzero}
\end{equation}
after the coordinate change~\eqref{change}.
\end{thm}

So, for each $\mu\in\cG$, the perturbed system~\eqref{perturbed} and the involution~\eqref{involution} admit a common analytic reducible invariant $n$-torus $\{\by=0, \, \bz=0\}$ with frequency vector $\tOmega(\mu)$ and Floquet matrix $\bfzero_m\oplus\tM(\mu)$, see~\eqref{withoutzero}. All such tori constitute an $s$-parameter Whitney $C^\infty$-smooth family.

\section{The Source Theorem in the Reversible Context~2}\label{source}

The material of this section almost coincides with that of Section~4 in \cite{S17MMJ}; we have included this section in the paper to achieve a self-contained presentation. To ``make'' the mapping~\eqref{mumapsto} submersive, one replaces $\Omega(\mu)+\Delta(\sigma,\mu)$ with an independent external parameter $\omega\in\mR^n$ and assumes the mapping
\[
\mu \mapsto \bigl(\alpha(\omega,\mu),\beta(\omega,\mu)\bigr)\in\mR^p
\]
(for $M$ dependent on $\omega$) to be submersive for fixed $\omega$.

Let $n\in\mZ_+$, $m\in\mN$, $p\in\mZ_+$, $s\in\mZ_+$, and $\omega_\star\in\mR^n$. Consider an analytic $(m+n+s)$-parameter family of analytic differential equations
\begin{equation}
\begin{aligned}
\dot{x} &= \omega+\xi(y,z,\sigma,\omega,\mu)+f(x,y,z,\sigma,\omega,\mu), \\
\dot{y} &= \sigma+\eta(y,z,\sigma,\omega,\mu)+g(x,y,z,\sigma,\omega,\mu), \\
\dot{z} &= M(\omega,\mu)z+\zeta(y,z,\sigma,\omega,\mu)+h(x,y,z,\sigma,\omega,\mu),
\end{aligned}
\label{BHTlike}
\end{equation}
where $x\in\mT^n$, $y\in\cO_m(0)$, $z\in\cO_{2p}(0)$ are the phase space variables, $\sigma\in\cO_m(0)$, $\omega\in\cO_n(\omega_\star)$, $\mu\in\cO_s(0)$ are external parameters, $M$ is a $2p\times 2p$ matrix-valued function, and $\xi=O(y,z)$, $\eta=O_2(y,z)$, $\zeta=O_2(y,z,\sigma)$. The functions $M$, $\xi$, $\eta$, $\zeta$ are supposed to be fixed whereas the terms $f$, $g$, $h$ are small perturbations. Let the systems~\eqref{BHTlike} be reversible with respect to the phase space involution~\eqref{involution}, where $R\in\GL(2p,\mR)$ is an involutive matrix with eigenvalues $1$ and $-1$ of multiplicity $p$ each, $M(\omega,\mu)R\equiv-RM(\omega,\mu)$, and the spectrum of $M(\omega_\star,0)$ is simple. One may assume that the spectrum of $M(\omega,\mu)$ is simple for any $\omega$ and $\mu$ and has the form $\fM\bigl(\nu_1,\nu_2,\nu_3;\alpha(\omega,\mu),\beta(\omega,\mu)\bigr)$ where $\nu_1+\nu_2+2\nu_3=p$ (see Definition~\ref{fM}). Retain the notation $\nu=\nu_2+\nu_3\in\mZ_+$.

\begin{thm}[\cite{S16RCD}]\label{thsource}
Suppose that the mapping
\[
\mu \mapsto \bigl(\alpha(\omega_\star,\mu),\beta(\omega_\star,\mu)\bigr)\in\mR^p
\]
is \emph{submersive} at the origin $\mu=0$ \textup{(}so that $s\geq p$\textup{)}. Then there exists a neighborhood $\fO\subset\mR^{n+s}$ of the point $(\omega_\star,0)$ such that for any closed set $\Gamma\subset\fO$ that is diffeomorphic to an $(n+s)$-dimensional ball and contains the point $(\omega_\star,0)$ in its interior, the following holds. For every complex neighborhood
\[
\cC\subset(\mC/2\pi\mZ)^n\times\mC^{2m+2p+n+s}
\]
of the set
\[
\mT^n\times\{0\in\mR^m\}\times\{0\in\mR^{2p}\}\times\{0\in\mR^m\}\times\Gamma
\]
and every $\cL\in\mN$, $\vare>0$, $\tau>n-1$ \textup{(}$\tau\geq 0$ for $n=0$\textup{)}, $\gamma>0$, there is a number $\delta>0$ with the following properties.

Suppose that the perturbation terms $f$, $g$, $h$ in~\eqref{BHTlike} can be holomorphically continued to the neighborhood $\cC$ and $|f|<\delta$, $|g|<\delta$, $|h|<\delta$ in $\cC$. Then for each $(\omega_0,\mu_0)\in\Gamma$, there exist points
\begin{equation}
v(\omega_0,\mu_0)\in\mR^m, \quad u(\omega_0,\mu_0)\in\mR^n, \quad w(\omega_0,\mu_0)\in\mR^s
\label{sdvig}
\end{equation}
and a change of variables
\begin{equation}
\begin{aligned}
x &= \bx+X(\bx,\omega_0,\mu_0), \\
y &= \by+Y^0(\bx,\omega_0,\mu_0)+Y^1(\bx,\omega_0,\mu_0)\by+Y^2(\bx,\omega_0,\mu_0)\bz, \\
z &= \bz+Z^0(\bx,\omega_0,\mu_0)+Z^1(\bx,\omega_0,\mu_0)\by+Z^2(\bx,\omega_0,\mu_0)\bz
\end{aligned}
\label{transformation}
\end{equation}
with $\bx\in\mT^n$, $\by\in\cO_m(0)$, $\bz\in\cO_{2p}(0)$ such that the following is valid.

\textup{1)} The functions $u$, $v$, $w$ in~\eqref{sdvig} are $C^\infty$-smooth as functions in $(\omega_0,\mu_0)$, and all the partial derivatives of each component of these functions of any order from $0$ to $\cL$ are smaller than $\vare$ in absolute value everywhere in $\Gamma$. The coefficients $X$, $Y^0$, $Y^1$, $Y^2$, $Z^0$, $Z^1$, $Z^2$ in~\eqref{transformation} are mappings ranging in $\mR^n$, $\mR^m$, $\gl(m,\mR)$, $\mR^{m\times 2p}$, $\mR^{2p}$, $\mR^{2p\times m}$, $\gl(2p,\mR)$, respectively. These mappings are analytic in $\bx$ and $C^\infty$-smooth in $(\omega_0,\mu_0)$. All the partial derivatives of each component of these mappings of any order from $0$ to $\cL$ are smaller than $\vare$ in absolute value everywhere in $\mT^n\times\Gamma$.

\textup{2)} For each $(\omega_0,\mu_0)\in\Gamma$, the change of variables~\eqref{transformation} commutes with the involution~\eqref{involution}.

\textup{3)} For \emph{any} point $(\omega_0,\mu_0)\in\Gamma$ such that the pair of vectors $\omega_0\in\mR^n$, $\beta(\omega_0,\mu_0)\in\mR^\nu$ is affinely $(\tau,\gamma,2)$-Diophantine \textup{(}see Definition~\ref{Diophantine}\textup{)}, the system~\eqref{BHTlike} at the parameter values
\begin{equation}
\sigma=v(\omega_0,\mu_0), \quad \omega=\omega_0+u(\omega_0,\mu_0), \quad \mu=\mu_0+w(\omega_0,\mu_0)
\label{shift}
\end{equation}
takes the form
\begin{equation}
\dot{\bx}=\omega_0+O(\by,\bz), \quad \dot{\by}=O_2(\by,\bz), \quad \dot{\bz}=M(\omega_0,\mu_0)\bz+O_2(\by,\bz)
\label{ideal}
\end{equation}
after the coordinate transformation~\eqref{transformation}.
\end{thm}

In fact, Theorem~\ref{thsource} is just a particular case of the main result of \cite{S16RCD}, see a discussion in \cite{S17MMJ}. Consider any point $(\omega_0,\mu_0)\in\Gamma$ such that the pair of vectors $\omega_0\in\mR^n$, $\beta(\omega_0,\mu_0)\in\mR^\nu$ is affinely $(\tau,\gamma,2)$-Diophantine. The perturbed system~\eqref{BHTlike} at the \emph{shifted} parameter values~\eqref{shift} has the reducible invariant $n$-torus $\{\by=0, \, \bz=0\}$ with \emph{the same} frequency vector $\omega_0$ and Floquet matrix $\bfzero_m\oplus M(\omega_0,\mu_0)$, see~\eqref{ideal}, as those of the reducible invariant $n$-torus $\{y=0, \, z=0\}$ of the system~\eqref{BHTlike} without the terms $f$, $g$, $h$ (the unperturbed system) at the parameter values $\sigma=0$, $\omega=\omega_0$, $\mu=\mu_0$. The torus $\{\by=0, \, \bz=0\}$ is analytic and invariant under the involution~\eqref{involution} and depends on $(\omega_0,\mu_0)$ in a $C^\infty$-way in the sense of Whitney.

\section{A Proof of Theorem~\ref{thmain}}\label{proof}

Our goal is to deduce Theorem~\ref{thmain} from Theorem~\ref{thsource} following the general Herman-like scheme (see \cite{S07TMIS} for a similar reduction technique in the ``conventional'' KAM contexts). Let the systems~\eqref{perturbed} satisfy the hypotheses of Theorem~\ref{thmain}. Since $M(\mu)$ depends on $\mu$ analytically and the spectrum of $M(0)$ is simple, one can introduce an additional parameter $\chi\in\cO_S(0)$ for an appropriate $S\in\mZ_+$ and construct an analytic family $M\new(\mu,\chi)$ of $2p\times 2p$ real matrices such that the following holds.

(1) $M\new(\mu,0)\equiv M(\mu)$ and $M\new(\mu,\chi)R\equiv-RM\new(\mu,\chi)$. As a consequence, one may assume that for any $\mu$ and $\chi$, the spectrum of $M\new(\mu,\chi)$ is simple and has the form
\[
\fM\bigl(\nu_1,\nu_2,\nu_3;\alpha\new(\mu,\chi),\beta\new(\mu,\chi)\bigr),
\]
where $\alpha\new(\mu,0)\equiv\alpha(\mu)$ and $\beta\new(\mu,0)\equiv\beta(\mu)$.

(2) The mapping
\[
(\mu,\chi) \mapsto \bigl(\alpha\new(\mu,\chi),\beta\new(\mu,\chi)\bigr)\in\mR^p
\]
is \emph{submersive} at $\mu=0$, $\chi=0$ (so that $s+S\geq p$).

The existence of a $2p\times 2p$ matrix-valued function $M\new$ satisfying these conditions follows immediately from the theory of normal forms and versal unfoldings of infinitesimally reversible matrices \cite{H96JDE,S91TSP,S93CJM}. It always suffices to set $S=p$.

Now introduce one more additional parameter $\omega\in\cO_n\bigl(\Omega(0)\bigr)$ and consider the analytic $(m+s+S+n)$-parameter family of analytic differential equations
\begin{equation}
\begin{aligned}
\dot{x} &= \omega+\xi(y,z,\sigma,\mu)+f(x,y,z,\sigma,\mu), \\
\dot{y} &= \sigma+\eta(y,z,\sigma,\mu)+g(x,y,z,\sigma,\mu), \\
\dot{z} &= M\new(\mu,\chi)z+\zeta(y,z,\sigma,\mu)+h(x,y,z,\sigma,\mu).
\end{aligned}
\label{extended}
\end{equation}
The systems~\eqref{extended} are reversible with respect to the involution~\eqref{involution} and satisfy all the hypotheses of Theorem~\ref{thsource}, with $\Omega(0)$, $s+S$, $(\mu,\chi)$, $M\new$ playing the roles of $\omega_\star$, $s$, $\mu$, $M$, respectively.

Consider a closed ball $A\subset\mR^{s-d}$ centered at the origin, a closed ball $B\subset\mR^d$ centered at the point $\fP_0$, a closed ball $\Gamma_1\subset\mR^S$ centered at the origin, and a closed ball $\Gamma_2\subset\mR^n$ centered at the point $\Omega(0)$. If the balls $A$ and $B$ are small enough then the set $\Gamma$~\eqref{Gamma} is well defined (and diffeomorphic to a closed $s$-dimensional ball). According to Theorem~\ref{thsource}, if all the four balls $A$, $B$, $\Gamma_1$, $\Gamma_2$ are sufficiently small then for every complex neighborhood~\eqref{neighborhood} of the set~\eqref{ostov} and every $\cL\in\mN$, $\tau>n-1$ ($\tau\geq 0$ for $n=0$), $\gamma>0$, the following holds.

Suppose that the perturbation terms $f$, $g$, $h$ in~\eqref{perturbed} and~\eqref{extended} can be holomorphically continued to the neighborhood~\eqref{neighborhood} and are sufficiently small in~\eqref{neighborhood}. Then for any $\mu_0\in\Gamma$, $\chi_0\in\Gamma_1$, and $\omega_0\in\Gamma_2$, there exist points
\begin{equation}
\begin{gathered}
v(\omega_0,\mu_0,\chi_0)\in\mR^m, \quad u(\omega_0,\mu_0,\chi_0)\in\mR^n, \\
w(\omega_0,\mu_0,\chi_0)\in\mR^s, \quad W(\omega_0,\mu_0,\chi_0)\in\mR^S
\end{gathered}
\label{bigsdvig}
\end{equation}
and a change of variables
\begin{equation}
\begin{aligned}
x &= \bx+\fX(\bx,\omega_0,\mu_0,\chi_0), \\
y &= \by+\fY^0(\bx,\omega_0,\mu_0,\chi_0)+\fY^1(\bx,\omega_0,\mu_0,\chi_0)\by+\fY^2(\bx,\omega_0,\mu_0,\chi_0)\bz, \\
z &= \bz+\fZ^0(\bx,\omega_0,\mu_0,\chi_0)+\fZ^1(\bx,\omega_0,\mu_0,\chi_0)\by+\fZ^2(\bx,\omega_0,\mu_0,\chi_0)\bz
\end{aligned}
\label{bigtransformation}
\end{equation}
with $\bx\in\mT^n$, $\by\in\cO_m(0)$, $\bz\in\cO_{2p}(0)$ such that the following is valid.

First, the functions $u$, $v$, $w$, $W$ in~\eqref{bigsdvig} are $C^\infty$-smooth. The coefficients $\fX$, $\fY^0$, $\fY^1$, $\fY^2$, $\fZ^0$, $\fZ^1$, $\fZ^2$ in~\eqref{bigtransformation} are analytic in $\bx$ and $C^\infty$-smooth in $(\omega_0,\mu_0,\chi_0)$. All the mappings $u$, $v$, $w$, $W$, $\fX$, $\fY^0$, $\fY^1$, $\fY^2$, $\fZ^0$, $\fZ^1$, $\fZ^2$ are small in the $C^{\cL}$-topology.

Second, for any $\mu_0\in\Gamma$, $\chi_0\in\Gamma_1$, and $\omega_0\in\Gamma_2$, the change of variables~\eqref{bigtransformation} commutes with the involution~\eqref{involution}.

Third, for \emph{any} points $\mu_0\in\Gamma$, $\chi_0\in\Gamma_1$, and $\omega_0\in\Gamma_2$ such that the pair of vectors $\omega_0\in\mR^n$, $\beta\new(\mu_0,\chi_0)\in\mR^\nu$ is affinely $(\tau,\gamma,2)$-Diophantine, the system~\eqref{extended} at the parameter values
\begin{equation}
\begin{gathered}
\sigma=v(\omega_0,\mu_0,\chi_0), \quad \omega=\omega_0+u(\omega_0,\mu_0,\chi_0), \\
\mu=\mu_0+w(\omega_0,\mu_0,\chi_0), \quad \chi=\chi_0+W(\omega_0,\mu_0,\chi_0)
\end{gathered}
\label{bigshift}
\end{equation}
takes the form
\begin{equation}
\dot{\bx}=\omega_0+O(\by,\bz), \quad \dot{\by}=O_2(\by,\bz), \quad \dot{\bz}=M\new(\mu_0,\chi_0)\bz+O_2(\by,\bz)
\label{bigideal}
\end{equation}
after the coordinate transformation~\eqref{bigtransformation}.

One may assume the balls $A$ and $B$ to be so small that $\Omega(\Gamma)$ lies in the interior of $\Gamma_2$. If the functions $u$, $v$, $w$, $W$ are small enough, then the system of equations
\begin{equation}
\begin{aligned}
\omega+u(\omega,\mu,\chi) &= \Omega\bigl(\mu+w(\omega,\mu,\chi)\bigr)+\Delta\bigl(v(\omega,\mu,\chi), \, \mu+w(\omega,\mu,\chi)\bigr), \\
\chi+W(\omega,\mu,\chi) &= 0
\end{aligned}
\label{keysystem}
\end{equation}
with $\mu\in\Gamma$ can be solved with respect to $\omega$ and $\chi$:
\[
\omega=\varphi(\mu), \quad \chi=\psi(\mu),
\]
where $\varphi:\Gamma\to\Gamma_2$ and $\psi:\Gamma\to\Gamma_1$ are $C^\infty$-functions close to $\Omega$ and $0$, respectively, in the $C^{\cL}$-topology. The key observation is that for any $\mu_0\in\Gamma$, the system~\eqref{extended} at the parameter values~\eqref{bigshift} with $\omega_0=\varphi(\mu_0)$ and $\chi_0=\psi(\mu_0)$ coincides with the \emph{original} system~\eqref{perturbed} at the parameter values
\[
\sigma=v(\omega_0,\mu_0,\chi_0), \quad \mu=\mu_0+w(\omega_0,\mu_0,\chi_0).
\]
Indeed, if $\omega_0=\varphi(\mu_0)$ and $\chi_0=\psi(\mu_0)$ then the equations~\eqref{keysystem} imply that the values of the parameters $\sigma$, $\omega$, $\mu$, $\chi$ given by~\eqref{bigshift} satisfy the relations
\[
\omega=\Omega(\mu)+\Delta(\sigma,\mu), \quad \chi=0.
\]

Let $\vare_3\in(0,1)$. Consider the closed $(s-d)$-dimensional ball $A'\subset A$ centered at the origin and the closed $d$-dimensional ball $B'\subset B$ centered at the point $\fP_0$ such that
\[
\meas_{s-d}A' = (1-\vare_3)^{1/2}\meas_{s-d}A, \quad \meas_dB' = (1-\vare_3)^{1/2}\meas_dB
\]
and set
\[
\Gamma'=\bigl\{ \mu(a,b) \bigm| a\in A', \, b\in B' \bigr\}\subset\Gamma
\]
($0\in\Gamma'$). If the functions $u$, $v$, $w$, $W$ are small enough, then the equation
\[
\mu=\mu_0+w\bigl(\varphi(\mu_0),\mu_0,\psi(\mu_0)\bigr)
\]
with $\mu\in\Gamma'$ can be solved with respect to $\mu_0$:
\[
\mu_0=\Upsilon(\mu),
\]
where $\Upsilon:\Gamma'\to\Gamma$ is a $C^\infty$-function close to the identity mapping $\mu\mapsto\mu$ in the $C^{\cL}$-topology.

We have arrived at the following conclusion. For any point $\mu\in\Gamma'$, set
\[
\mu_0=\Upsilon(\mu), \quad \omega_0=\varphi\bigl(\Upsilon(\mu)\bigr), \quad \chi_0=\psi\bigl(\Upsilon(\mu)\bigr).
\]
If the pair of vectors $\omega_0$, $\beta\new(\mu_0,\chi_0)$ is affinely $(\tau,\gamma,2)$-Diophantine, then the \emph{original} system~\eqref{perturbed} at the parameter values $\mu$ and $\sigma=v(\omega_0,\mu_0,\chi_0)$ takes the form~\eqref{bigideal} after the coordinate transformation~\eqref{bigtransformation}.

Introduce the functions
\begin{align*}
\hOmega(\mu) &= \varphi\bigl(\Upsilon(\mu)\bigr), \\
\Psi(\mu) &= \psi\bigl(\Upsilon(\mu)\bigr), \\
\hM(\mu) &= M\new\bigl(\Upsilon(\mu),\Psi(\mu)\bigr), \\
\halpha(\mu) &= \alpha\new\bigl(\Upsilon(\mu),\Psi(\mu)\bigr), \\
\hbeta(\mu) &= \beta\new\bigl(\Upsilon(\mu),\Psi(\mu)\bigr), \\
\hTheta(\mu) &= v\Bigl(\hOmega(\mu),\Upsilon(\mu),\Psi(\mu)\Bigr)
\end{align*}
for $\mu\in\Gamma'$ and
\begin{align*}
\hX(\bx,\mu) &= \fX\Bigl(\bx,\hOmega(\mu),\Upsilon(\mu),\Psi(\mu)\Bigr), \\
\hY^r(\bx,\mu) &= \fY^r\Bigl(\bx,\hOmega(\mu),\Upsilon(\mu),\Psi(\mu)\Bigr), \quad r=0,1,2, \\
\hZ^r(\bx,\mu) &= \fZ^r\Bigl(\bx,\hOmega(\mu),\Upsilon(\mu),\Psi(\mu)\Bigr), \quad r=0,1,2
\end{align*}
for $\bx\in\mT^n$ and $\mu\in\Gamma'$. The mappings $\hOmega$, $\Psi$, $\hM$, $\halpha$, $\hbeta$, $\hTheta$ are $C^\infty$-smooth and the functions $\hOmega-\Omega$, $\Psi$, $\hM-M$, $\halpha-\alpha$, $\hbeta-\beta$, $\hTheta$ are small in the $C^{\cL}$-topology. For any $\mu\in\Gamma'$, one has $\hM(\mu)R=-R\hM(\mu)$, and the spectrum of the $2p\times 2p$ matrix $\hM(\mu)$ is simple and has the form $\fM\Bigl(\nu_1,\nu_2,\nu_3;\halpha(\mu),\hbeta(\mu)\Bigr)$. The coefficients $\hX$, $\hY^0$, $\hY^1$, $\hY^2$, $\hZ^0$, $\hZ^1$, $\hZ^2$ are analytic in $\bx\in\mT^n$, $C^\infty$-smooth in $\mu\in\Gamma'$, and small in the $C^{\cL}$-topology provided that the perturbation terms $f$, $g$, $h$ in~\eqref{perturbed} are small enough.

The conclusion we have come to so far can be reformulated as follows. If the pair of vectors $\hOmega(\mu)\in\mR^n$, $\hbeta(\mu)\in\mR^\nu$ is affinely $(\tau,\gamma,2)$-Diophantine for some $\mu\in\Gamma'$, then the system~\eqref{perturbed} at the parameter values $\mu$ and $\sigma=\hTheta(\mu)$ takes the form
\[
\dot{\bx}=\hOmega(\mu)+O(\by,\bz), \quad \dot{\by}=O_2(\by,\bz), \quad \dot{\bz}=\hM(\mu)\bz+O_2(\by,\bz)
\]
after the $G$-commuting coordinate change
\begin{align*}
x &= \bx+\hX(\bx,\mu), \\
y &= \by+\hY^0(\bx,\mu)+\hY^1(\bx,\mu)\by+\hY^2(\bx,\mu)\bz, \\
z &= \bz+\hZ^0(\bx,\mu)+\hZ^1(\bx,\mu)\by+\hZ^2(\bx,\mu)\bz
\end{align*}
with $\bx\in\mT^n$, $\by\in\cO_m(0)$, $\bz\in\cO_{2p}(0)$.

Consider the closed $(s-d)$-dimensional balls $\tA\subset A''\subset A'$ centered at the origin and the closed $d$-dimensional ball $\tB\subset B'$ centered at the point $\fP_0$ such that
\[
\meas_{s-d}A'' = (1-\vare_3)^{3/4}\meas_{s-d}A
\]
and the relations~\eqref{tAtB} hold. Define the sets $\tGamma\subset\Gamma''\subset\Gamma'$ by the equation
\[
\Gamma''=\Bigl\{ \mu(a,b) \Bigm| a\in A'', \, b\in\tB \Bigr\}
\]
($0\in\Gamma''$) and the equation~\eqref{tGamma}. Let the balls $A$ and $B$ be so small that the Jacobian~\eqref{Jacobian} vanishes nowhere in $\Gamma$ for $d\geq 1$. Then the system of equations
\begin{align*}
\hOmega_+(\mu^\ast_+,\mu_-) &= \Omega_+(\mu_+,\mu_-), \\
\halpha_+(\mu^\ast_+,\mu_-) &= \alpha_+(\mu_+,\mu_-), \\
\hbeta_+(\mu^\ast_+,\mu_-) &= \beta_+(\mu_+,\mu_-)
\end{align*}
can be solved with respect to $\mu^\ast_+$ for $\mu=(\mu_+,\mu_-)\in\Gamma''$ provided that the functions $u$, $v$, $w$, $W$ are small enough:
\[
\mu^\ast_+=\mu_++\Xi_+(\mu_+,\mu_-)=\mu_++\Xi_+(\mu),
\]
where $\Xi_+:\Gamma''\to\mR^d$ is a $C^\infty$-function small in the $C^{\cL}$-topology. ``Complement'' the mapping $\Xi_+$ with the \emph{zero} function $\Xi_-:\Gamma''\to\mR^{s-d}$ in such a way that
\[
\Xi_+=(\Xi_l \mid l\in\fT), \quad \Xi_-=(\Xi_l \mid l\notin\fT)
\]
for the mapping $\Xi=(\Xi_+,\Xi_-):\Gamma''\to\mR^s$, cf.~\eqref{fT}. If $d=0$ then $\Xi\equiv 0$. For $\mu\in\Gamma''$ one has $\mu+\Xi(\mu)\in\Gamma'$ and
\begin{align*}
\hOmega_+\bigl(\mu+\Xi(\mu)\bigr) &= \Omega_+(\mu), \\
\halpha_+\bigl(\mu+\Xi(\mu)\bigr) &= \alpha_+(\mu), \\
\hbeta_+\bigl(\mu+\Xi(\mu)\bigr) &= \beta_+(\mu).
\end{align*}

Now set
\begin{align*}
\tOmega(\mu) &= \hOmega\bigl(\mu+\Xi(\mu)\bigr), \\
\tM(\mu) &= \hM\bigl(\mu+\Xi(\mu)\bigr), \\
\talpha(\mu) &= \halpha\bigl(\mu+\Xi(\mu)\bigr), \\
\tbeta(\mu) &= \hbeta\bigl(\mu+\Xi(\mu)\bigr), \\
\Theta(\mu) &= \hTheta\bigl(\mu+\Xi(\mu)\bigr)
\end{align*}
for $\mu\in\Gamma''$ and
\begin{align*}
X(\bx,\mu) &= \hX\bigl(\bx,\mu+\Xi(\mu)\bigr), \\
Y^r(\bx,\mu) &= \hY^r\bigl(\bx,\mu+\Xi(\mu)\bigr), \quad r=0,1,2, \\
Z^r(\bx,\mu) &= \hZ^r\bigl(\bx,\mu+\Xi(\mu)\bigr), \quad r=0,1,2
\end{align*}
for $\bx\in\mT^n$ and $\mu\in\Gamma''$. The mappings $\Xi$, $\tOmega$, $\tM$, $\talpha$, $\tbeta$, $\Theta$ are $C^\infty$-smooth and the functions $\Xi$, $\tOmega-\Omega$, $\tM-M$, $\talpha-\alpha$, $\tbeta-\beta$, $\Theta$ are small in the $C^{\cL}$-topology. For any $\mu\in\Gamma''$, one has $\tM(\mu)R=-R\tM(\mu)$, and the spectrum of the $2p\times 2p$ matrix $\tM(\mu)$ is simple and has the form $\fM\Bigl(\nu_1,\nu_2,\nu_3;\talpha(\mu),\tbeta(\mu)\Bigr)$. The identities~\eqref{identities} are valid in $\Gamma''$. The coefficients $X$, $Y^0$, $Y^1$, $Y^2$, $Z^0$, $Z^1$, $Z^2$ are analytic in $\bx\in\mT^n$, $C^\infty$-smooth in $\mu\in\Gamma''$, and small in the $C^{\cL}$-topology provided that the perturbation terms $f$, $g$, $h$ in~\eqref{perturbed} are small enough. For each $\mu\in\Gamma''$, the change of variables~\eqref{change} commutes with the involution~\eqref{involution}.

For any $\mu^0\in\Gamma''$ such that the pair of vectors $\tOmega(\mu^0)$, $\tbeta(\mu^0)$ is affinely $(\tau,\gamma,2)$-Diophantine, the system~\eqref{perturbed} at the parameter values $\mu=\mu^0+\Xi(\mu^0)$ and $\sigma=\Theta(\mu^0)$ takes the form~\eqref{goal} after the coordinate change~\eqref{change} with $\mu=\mu^0$.

The pair of mappings~\eqref{pair} is affinely $(Q,2)$-nondegenerate at $a=0$. One may assume the balls $A$ and $B$ (and, consequently, the balls $\tA$ and $\tB$) to be so small that the pair of mappings
\begin{equation}
a\mapsto\Omega_-\bigl(\mu(a,b)\bigr)\in\mR^{n-d_1} \quad \text{and} \quad a\mapsto\beta_-\bigl(\mu(a,b)\bigr)\in\mR^{\nu-d_3}
\label{ourpair}
\end{equation}
is affinely $(Q,2)$-nondegenerate at each point $a\in\tA$ for any fixed value of $b\in\tB$. Now the Diophantine Lemma~\ref{lmDiophantine} can be applied with

$s-d\geq 1$, $n-d_1$, $\nu-d_3$, $d_1$, $d_3$, and $2$ playing the roles of $s$, $n$, $\nu$, $n\add$, $\nu\add$, and $L$, respectively,

$\tB$, $\tA$, and the interior of $A''$ playing the roles of $B$, $A$, and $\fA$, respectively,

the mappings~\eqref{ourpair} playing the roles of the mappings~\eqref{abstractpair},

the mappings $b\mapsto b^{:1}\in\mR^{d_1}$ and $b\mapsto b^{:3}\in\mR^{d_3}$ playing the roles of the mappings $\Omega\add$ and $\beta\add$, respectively,

$\vare_2$ playing the role of $\vare$.

According to Lemma~\ref{lmDiophantine}, if $\cL\geq Q$ and the differences $\tOmega-\Omega$ and $\tbeta-\beta$ are sufficiently small in $\Gamma''$ in the $C^{\cL}$-topology, then the following is valid. Let $\vare_2\in(0,1)$, $\tau_\ast\geq\max(0,d_1-1)$, and $\gamma_\ast>0$. Suppose that $\tau>nQ$, $\tau\geq\tau_\ast$, and $\gamma$ is sufficiently small: $0<\gamma\leq\gamma_0(\vare_2,\tau,\gamma_\ast)$. Then for any point $b\in\tB$ such that the pair of vectors $b^{:1}\in\mR^{d_1}$, $b^{:3}\in\mR^{d_3}$ is affinely $(\tau_\ast,\gamma_\ast,2)$-Diophantine, the Lebesgue measure $\meas_{s-d}$ of the set $\cG_b$ of those points $a\in\tA$ for which the pair of vectors $\tOmega\bigl(\mu(a,b)\bigr)\in\mR^n$, $\tbeta\bigl(\mu(a,b)\bigr)\in\mR^\nu$ is affinely $(\tau,\gamma,2)$-Diophantine, is greater than $(1-\vare_2)\meas_{s-d}\tA$. Indeed, denoting $\mu(a,b)$ by $\mu^0$ and taking into account the identities~\eqref{rasklad} and~\eqref{identities}, one has
\begin{gather*}
\tOmega(\mu^0) = \Bigl(\tOmega_+(\mu^0),\tOmega_-(\mu^0)\Bigr) = \Bigl(\Omega_+(\mu^0),\tOmega_-(\mu^0)\Bigr) = \Bigl(b^{:1},\tOmega_-(\mu^0)\Bigr), \\
\tbeta(\mu^0) = \Bigl(\tbeta_+(\mu^0),\tbeta_-(\mu^0)\Bigr) = \Bigl(\beta_+(\mu^0),\tbeta_-(\mu^0)\Bigr) = \Bigl(b^{:3},\tbeta_-(\mu^0)\Bigr).
\end{gather*}
The proof of Theorem~\ref{thmain} is completed.

\section{Invariant Tori in General Systems}\label{multiparameter}

Let $n\in\mZ_+$, $m\in\mN$, $p\in\mZ_+$. Consider a system of differential equations
\begin{equation}
\dot{x}=\fx(x,y,z), \quad \dot{y}=\fy(y), \quad \dot{z}=\fz(x,y,z),
\label{general}
\end{equation}
where $x\in\mT^n$, $y\in\cO_m(0)$, $z\in\cO_{2p}(0)$ are the phase space variables, cf.~\eqref{unperturbed} and~\eqref{perturbed}, and suppose that this system is reversible with respect to the phase space involution~\eqref{involution}, where $R\in\GL(2p,\mR)$ is an involutive matrix. The reversibility condition means that
\[
\fx(-x,-y,Rz)\equiv\fx(x,y,z), \quad \fy(-y)\equiv\fy(y), \quad \fz(-x,-y,Rz)\equiv-R\fz(x,y,z).
\]
Note that we impose no restrictions on the equations for $\dot{x}$ and $\dot{z}$ (apart from the reversibility) but the right-hand side of the equation for $\dot{y}$ is assumed to be independent of $x$ and $z$. In the present section, we give a rigorous proof of the following statement.

\begin{prp}\label{notorus}
Let the system~\eqref{general} and the involution~\eqref{involution} admit a common invariant torus carrying quasi-periodic motions. Then $\fy(0)=0$.
\end{prp}

The torus in Proposition~\ref{notorus} is not assumed to be of dimension $n$ (not to mention to be close to the torus $\{y=0, \, z=0\}$).

In particular, suppose that the system~\eqref{general} depends on a $\fc$-dimensional parameter $\fw$ with $\fc<m$:
\[
\dot{x}=\fx(x,y,z,\fw), \quad \dot{y}=\fy(y,\fw), \quad \dot{z}=\fz(x,y,z,\fw).
\]
Generically the points $\fy(0,\fw)$ constitute a $\fc$-dimensional surface in $\mR^m$ which does not contain the origin. Consequently, if $\fc<m$ then a generic $\fc$-parameter family of $G$-reversible systems~\eqref{general} admits an invariant torus carrying quasi-periodic motions (and invariant under the involution $G$ as well) at \emph{no} value of the parameter.

Proposition~\ref{notorus} is a particular case of the following more general statement.

\begin{prp}\label{notori}
Let the system
\begin{equation}
\dot{u}=U(u,v), \qquad \dot{v}=V(v)
\label{system}
\end{equation}
of differential equations on the direct product $A\times B=\bigl\{(u,v)\bigr\}$ of manifolds $A$ and $B$ be reversible with respect to an involution $G:(u,v)\mapsto\bigl(G_A(u),G_B(v)\bigr)$ where $G_A$ and $G_B$ are involutions of $A$ and $B$, respectively. Suppose that $\Fix G_B$ consists of a single point $v^0\in B$. Let the system~\eqref{system} and the involution $G$ admit a common invariant torus carrying quasi-periodic motions. Then $V(v^0)=0$.
\end{prp}

In turn, the proof of Proposition~\ref{notori} is based on the following lemma.

\begin{lmm}\label{fundamental}
Let $F:\mT^n\to K$ be a surjective continuous mapping of $\mT^n$ onto a compact topological space. Let $\fg^t$ be a quasi-periodic flow on $\mT^n$, i.e., $\fg^t(\phi)=\phi+\omega t$ \textup{(}$\phi\in\mT^n$\textup{)} where $\omega\in\mR^n$ is a fixed vector with rationally independent components. Let also $\fG^t$ be a continuous action of $\mR$ on $K$. Suppose that $F\circ\fg^t=\fG^t\circ F$. Then $K$ is a torus of dimension no greater than $n$, and $\fG^t$ is quasi-periodic.
\end{lmm}

It is hardy possible that Lemma~\ref{fundamental} is new, but I have failed to find it in the literature.

Let us first deduce Proposition~\ref{notori} from Lemma~\ref{fundamental}. Suppose that the system~\eqref{system} and the involution $G$ admit a common invariant $n$-torus $F(\mT^n)$, where $F=(F_A,F_B)$ is an embedding of $\mT^n$ into $A\times B$ ($F_A$ and $F_B$ take $\mT^n$ to $A$ and $B$, respectively). Assume the torus $F(\mT^n)$ to carry a quasi-periodic flow $F\circ\fg^t\circ F^{-1}$ where $\fg^t$ is a quasi-periodic flow on $\mT^n$. Since $F(\mT^n)$ is invariant under $G$, for each $\phi\in\mT^n$ there exists $\phi'\in\mT^n$ such that $F_A(\phi')=G_A\bigl(F_A(\phi)\bigr)$ and $F_B(\phi')=G_B\bigl(F_B(\phi)\bigr)$. Consequently, the sets $F_A(\mT^n)$ and $F_B(\mT^n)$ are invariant under the involutions $G_A$ and $G_B$, respectively. It is also clear that $F_B(\mT^n)$ is an invariant set of the equation $\dot{v}=V(v)$ and $F_B\circ\fg^t=\fG^t\circ F_B$ where $\fG^t$ is the restriction of the flow of the vector field $V$ to $F_B(\mT^n)$. According to Lemma~\ref{fundamental}, $F_B(\mT^n)$ is a $k$-torus for some $k$ ($0\leq k\leq n$), and $\fG^t$ is a quasi-periodic flow on $F_B(\mT^n)$. According to Lemma~\ref{torus}, the torus $F_B(\mT^n)$ contains $2^k$ fixed points of the involution $G_B$. Since $\Fix G_B=\{v^0\}$, we arrive at the conclusion that $k=0$ and $F_B(\mT^n)=\{v^0\}$. Now the invariance of the set $F_B(\mT^n)$ under the flow of $V$ implies that $V(v^0)=0$.

It remains to prove Lemma~\ref{fundamental}. Let $\Lambda=F^{-1}\bigl(F(0)\bigr)\subset\mT^n$. Then $\Lambda$ is a closed subset of $\mT^n$. Our first goal is to verify that $F(\phi^1)=F(\phi^2)$ if and only if $\phi^1-\phi^2\in\Lambda$. Indeed, consider the sequence $\{t_j\}_{j\in\mN}$ of real numbers such that $\lim_{j\to\infty} \omega t_j=\phi^2$ on $\mT^n$. If $F(\phi^1)=F(\phi^2)$, then
\begin{multline*}
F(\phi^1-\phi^2)=\lim_{j\to\infty} F(\phi^1-\omega t_j)=\lim_{j\to\infty} \fG^{-t_j}\bigl(F(\phi^1)\bigr) \\
{}=\lim_{j\to\infty} \fG^{-t_j}\bigl(F(\phi^2)\bigr)=\lim_{j\to\infty} F(\phi^2-\omega t_j)=F(0),
\end{multline*}
so that $\phi^1-\phi^2\in\Lambda$. On the other hand, if $\phi^1-\phi^2\in\Lambda$, i.e., $F(\phi^1-\phi^2)=F(0)$, then
\begin{multline*}
F(\phi^1)=\lim_{j\to\infty} F(\phi^1-\phi^2+\omega t_j)=\lim_{j\to\infty} \fG^{t_j}\bigl(F(\phi^1-\phi^2)\bigr) \\
{}=\lim_{j\to\infty} \fG^{t_j}\bigl(F(0)\bigr)=\lim_{j\to\infty} F(\omega t_j)=F(\phi^2).
\end{multline*}
In particular, if $\phi^1\in\Lambda$ and $\phi^2\in\Lambda$, then $F(\phi^1+\phi^2)=F(\phi^1)=F(0)$, so that $\phi^1+\phi^2\in\Lambda$, and $F(-\phi^1)=F(0)$, so that $-\phi^1\in\Lambda$. Thus, $\Lambda$ is a closed \emph{subgroup} of $\mT^n$, and there exists a natural bijection $\mT^n/\Lambda\to F(\mT^n)$.

Now one can apply the Pontryagin duality theorem for locally compact Abelian groups (also known as the Pontryagin--van Kampen duality theorem). For the consequences of this fundamental theorem we need and for their particular case that concerns $\mT^n$, the reader is referred to, e.g., Proposition~38 in \cite{M77} and Corollary~1.2.2 on page~706 in \cite{T13}. According to the Pontryagin duality theorem, closed subgroups of $\mT^n$ are characterized by their \emph{annihilators} in the character group $X^\ast(\mT^n)\approx\mZ^n$, i.e., in the group of all the continuous homomorphisms $\chi:\mT^n\to\mS^1=\mR/2\pi\mZ$:
\[
\chi(\phi_1,\ldots,\phi_n)=m_1\phi_1+\cdots+m_n\phi_n, \qquad m_1,\ldots,m_n\in\mZ
\]
(the groups $\mT^n$ and $\mZ^n$ are \emph{dual} to each other). In other words, there exists a subgroup $L\subset\mZ^n$ such that
\[
\Lambda=\bigl\{ \phi\in\mT^n \bigm|
m_1\phi_1+\cdots+m_n\phi_n=0 \;\; \forall(m_1,\ldots,m_n)\in L \bigr\}.
\]

On the other hand, there is a matrix $Q\in\SL(n,\mZ)$ such that
\[
LQ=\bigl\{ (\,\underbrace{0,\ldots,0}_{n-k}\,, \, q_1r_1,\ldots,q_kr_k) \bigm|
r_1,\ldots,r_k\in\mZ \bigr\},
\]
where $q_1\geq\cdots\geq q_k\geq 1$ are certain natural numbers (here the elements of $\mZ^n$ are regarded as row vectors and $0\leq k\leq n$). The rank of the lattice $L$ is equal to $k$. Introduce the new coordinate frame $\psi=Q^{-1}\phi$ on the torus $\mT^n$ (here the points of $\mT^n$ are regarded as column vectors). In the new coordinate frame
\[
\Lambda=\bigl\{ (\psi_1,\ldots,\psi_{n-k}, \, 2\pi p_1/q_1,\ldots,2\pi p_k/q_k) \bigr\},
\]
where
\[
0\leq p_1\leq q_1-1,\ldots,0\leq p_k\leq q_k-1; \qquad
\psi_1,\ldots,\psi_{n-k}\in\mS^1, \quad p_1,\ldots,p_k\in\mZ_+,
\]
and $\dim\Lambda=n-k$. Moreover, the flow $\fg^t$ on $\mT^n$ in the new coordinate frame is determined by the equation
\[
\dot{\psi}=Q^{-1}\dot{\phi}=Q^{-1}\omega=\varpi=(\varpi_1,\ldots,\varpi_n).
\]
The set $K=F(\mT^n)$ is homeomorphic to $\mT^n/\Lambda\approx\mT^k$. The natural coordinates on the factor $\mT^n/\Lambda$ and, consequently, on the set $F(\mT^n)$ are the coordinates
\[
(q_1\psi_{n-k+1},\ldots,q_k\psi_n)\in\mT^k.
\]
The flow $\fG^t$ on the $k$-torus $F(\mT^n)$ is quasi-periodic with the frequency vector
\[
(q_1\varpi_{n-k+1},\ldots,q_k\varpi_n)\in\mR^k.
\]
This completes the proof of Lemma~\ref{fundamental}.

\begin{rmk}\label{MathOverflow}
Lemma~\ref{fundamental} is probably related to the theory of minimal isometric systems, cf.\ Proposition~2.6.7 in the book \cite{T09}.
\end{rmk}

\end{document}